%% file: primary.tex
\documentclass[a4paper, 9pt, twosides]{amsart}

\usepackage{amssymb}
\usepackage{graphicx}
\usepackage{enumitem}
\usepackage[usenames, dvipsnames]{color}
\usepackage{todonotes}
\usepackage[colorlinks=true,
		raiselinks=true,
		linkcolor=MidnightBlue,
		citecolor=Mahogany,
		urlcolor=ForestGreen,
		pdfauthor={Yogesh Kumar, Sukumar Srikant and Debasish Chatterjee},
		pdftitle={Optimal multiplexing of sparse controllers},
		pdfkeywords={sparse control, nonsmooth Pontryagin maximum principle, networked control, optimal control},
		pdfsubject={Manuscript},
		plainpages=false]{hyperref}

\usepackage{newtxtext}
\usepackage[varg, bigdelims, vvarbb]{newtxmath}
\DeclareFontFamily{T1}{pzc}{}
\DeclareFontShape{T1}{pzc}{m}{it}{<-> [1.2] pzcmi8t}{}
\DeclareMathAlphabet{\mathpzc}{T1}{pzc}{m}{it}

\usepackage{multicol}
\usepackage{mathtools}
\mathtoolsset{mathic=true}

\allowdisplaybreaks
\newtheorem{theorem}{Theorem}[section]

\newtheorem{corollary}[theorem]{Corollary}

\theoremstyle{definition}

\newtheorem{example}[theorem]{Example}

\theoremstyle{remark}
\newtheorem{remark}[theorem]{Remark}

\theoremstyle{assumption}

\theoremstyle{fact}

\theoremstyle{claim}

\numberwithin{equation}{section}
\DeclareMathOperator*{\minimize}{minimize}
\DeclareMathOperator*{\argmax}{arg\,max}
\DeclareMathOperator*{\sbjto}{subject\ to}
\DeclareMathOperator{\Leb}{Leb}
\DeclareMathOperator{\sgn}{sgn}
\DeclareMathOperator{\blkdiag}{blkdiag}

\renewcommand{\mapsto}{\longmapsto}
\newcommand{\lra}{\longrightarrow}
\renewcommand{\le}{\leqslant}
\renewcommand{\ge}{\geqslant}
\newcommand{\Let}{\coloneqq}
\newcommand{\teL}{\eqqcolon}
\newcommand{\setmin}{\smallsetminus}
\newcommand{\ol}{\overline}
\newcommand{\wt}{\widetilde}

\newcommand{\R}{\mathbb{R}}

\newcommand{\admact}{\mathbb{U}}
\newcommand{\subadmact}[1]{\admact^{(#1)}}
\newcommand{\admactreach}{\admact_{\mathrm{rch}}}
\newcommand{\admactLQ}{\admact_{\mathrm{LQ}}}
\newcommand{\admactMayer}{\admact_{\mathrm{M}}}
\newcommand{\reachmux}{\sigma^{\mathrm{rch}}}
\newcommand{\reachmuxset}{\Sigma^{\mathrm{rch}}}
\newcommand{\LQmux}{\sigma^{\mathrm{LQ}}}
\newcommand{\LQmuxset}{\Sigma^{\mathrm{LQ}}}
\newcommand{\Mayermux}{\sigma^{\mathrm{M}}}

\newcommand{\dd}{\mathrm{d}}
\newcommand{\trnsp}{^\top}
\newcommand{\inverse}{^{-1}}
\newcommand{\opt}{_\ast}
\newcommand{\Ham}{H}

\newcommand{\tinit}{0}
\newcommand{\tfin}{\hat t}
\newcommand{\stinit}{\bar x}
\newcommand{\stfin}{\hat x}
\newcommand{\st}{x}
\newcommand{\subst}[1]{\st^{#1}}
\newcommand{\substinit}[1]{\bar{\st}^{#1}}
\newcommand{\substfin}[1]{\hat{\st}^{#1}}
\newcommand{\tdsubst}[1]{\dot{\st}^{#1}}
\newcommand{\con}{u}
\newcommand{\subcon}[1]{\con^{#1}}
\newcommand{\adj}{p}
\newcommand{\subadj}[1]{\adj^{#1}}
\newcommand{\tdsubadj}[1]{\dot{\adj}^{#1}}

\newcommand{\norm}[1]{\left\lVert{#1}\right\rVert}
\newcommand{\abs}[1]{\left\lvert{#1}\right\rvert}
\newcommand{\Lp}[1]{\mathpzc L_{#1}}
\newcommand{\indic}[1]{\mathsf{1}_{#1}}
\newcommand{\inprod}[2]{\left\langle{#1}, {#2}\right\rangle}

\newcommand{\pmat}[1]{\begin{pmatrix}#1\end{pmatrix}}

\newcommand{\secref}[1]{\S\ref{#1}}

	\title{Optimal multiplexing of sparse controllers for linear systems}
	\author[Y.\ Kumar]{Yogesh Kumar}
	\address{YK is with the Department of Aerospace Engineering, Indian Institute of Technology Bombay, Powai, Mumbai 400076, India.}
	\email{02yogesh16@gmail.com}
	\author[S.\ Srikant]{Sukumar Srikant}
	\author[D.\ Chatterjee]{Debasish Chatterjee}
	\address{SS and DC are with Systems \& Control Engineering, Indian Institute of Technology Bombay, Powai, Mumbai 400076, India.}
	\email{srikant.sukumar@iitb.ac.in, dchatter@iitb.ac.in}

\begin{document}

	%\date{\DTMnow}

	\begin{abstract}
		This article treats three problems of sparse and optimal multiplexing a finite ensemble of linear control systems. Given an ensemble of linear control systems, multiplexing of the controllers consists of an algorithm that selects, at each time \(t\), only one from the ensemble of linear systems is actively controlled whereas the other systems evolve in open-loop. The first problem treated here is a ballistic reachability problem where the control signals are required to be maximally sparse and multiplexed, the second concerns sparse and optimally multiplexed linear quadratic control, and the third is a sparse and optimally multiplexed Mayer problem. Numerical experiments are provided to demonstrate the efficacy of the techniques developed here.
	\end{abstract}

    \maketitle

	\begin{multicols}{2}

	\input{intro}

%%%%%%%%%%%%%%%%%%%%%%%%%%%%%%%%%%%%%%%%%%%%%%%%%%%%%%%%%%%%%%%%%%%%%%%%%%%%%%%%
	\section{Premise and preliminaries}
	\label{s:prelims}
%%%%%%%%%%%%%%%%%%%%%%%%%%%%%%%%%%%%%%%%%%%%%%%%%%%%%%%%%%%%%%%%%%%%%%%%%%%%%%%%
		Let \(N\) be a positive integer, and let \(\tfin > 0\) be a fixed time instant. Consider the finite ensemble of linear time-invariant control systems given by
		\begin{equation}
		\label{e:ensemble}
		\left\{
		\begin{aligned}
			& \tdsubst{k}(t) = A_k \subst{k}(t) + B_k \subcon{k}(t)\quad\text{for a.e.\ }t\in[\tinit, \tfin],\\
			& k = 1, \ldots, N,
		\end{aligned}
		\right.
		\end{equation}
		with the following data:
		\begin{enumerate}[label=(\ref{e:ensemble}-\alph*), leftmargin=*, widest=b, align=left]
			\item \label{e:ensemble:AB} \(A_k\in\R^{d_k\times d_k}\) and \(B_k\in\R^{d_k\times m_k}\),
			\item \label{e:ensemble:x} the \emph{states} \(\subst{k}(t)\in\R^{d_k}\) with \(\subst{k}(\tinit) = \substinit{k}\),
			\item \label{e:ensemble:admact} the \emph{control actions} \(\subcon{k}(t)\in\subadmact{k} \subset \R^{m_k}\), where the \emph{set of admissible control actions} \(\subadmact{k}\) contains \(0\in\R^{m_k}\), and
			\item \label{e:ensemble:controls} the \emph{controls} \([\tinit, \tfin]\ni t\mapsto \subcon{k}(t)\in\subadmact{k}\) are measurable.\footnote{For us the word `measurability' always refers to Lebesgue measurability, and `a.e.'\ will refer to almost everywhere relative to the Lebesgue measure.}
		\end{enumerate}
		Define \(d \Let \sum_{k=1}^N d_k\) and \(m \Let \sum_{k=1}^N m_k\). For a control \(\subcon{k}\) injected into the \(k\)-th control system in \eqref{e:ensemble}, we let \([\tinit, \tfin]\ni t\mapsto \subst{k}(t)\in\R^{d_k}\) denote the unique absolutely continuous solution of the \(k\)-th control system corresponding to \(\subcon{k}\); we call \(\subst{k}\) the \emph{state trajectory} corresponding to \(\subcon{k}\) (we suppress the explicit dependence of \(\subst{k}\) on \(\subcon{k}\) to reduce notational clutter).\footnote{For us the word `solution' always refers to a solution in the sense of Carath\'eodory \cite[Chapter 1]{ref:Fil-88}. Note that existence and uniqueness of Carath\'eodory solutions of each member of \eqref{e:ensemble} under measurable controls are guaranteed by linearity of the states on the right-hand side.} The map
		\[
			[\tinit, \tfin]\ni t\mapsto \bigl(\subst{k}(t), \subcon{k}(t)\bigr)\in\R^{d_k}\times\subadmact{k}
		\]
		is known as an \emph{admissible state-action trajectory} whenever \(\subst{k}\) is the solution of the \(k\)-th member of \eqref{e:ensemble} under \(\subcon{k}\). Out of the collection \(\bigl((\subst{k}, \subcon{k})\bigr)_{k=1}^N\) of state-action trajectories we define the two maps
		\begin{equation}
		\label{e:joint state-actions}
		\left\{
		\begin{aligned}
			& [\tinit, \tfin]\ni t\mapsto \st(t) \Let \pmat{\subst{1}(t)\\ \smash[t]{\vdots} \\ \subst{N}(t)} \in \R^d,\\
			& [\tinit, \tfin]\ni t\mapsto \con(t) \Let \pmat{\subcon{1}(t)\\ \smash[t]{\vdots} \\ \subcon{N}(t)} \in \subadmact{1}\times\cdots\times\subadmact{N}.
		\end{aligned}
		\right.
		\end{equation}
		With \(A \Let \blkdiag(A_1, \ldots, A_N)\) and \(B \Let \blkdiag(B_1, \ldots, B_N)\), the ensemble of systems \eqref{e:ensemble} admits a compact representation as the joint control system
		\begin{equation}
		\label{e:joint system}
			\dot\st(t) = A \st(t) + B \con(t)\quad \text{for a.e.\ }t\in[\tinit, \tfin].
		\end{equation}

		Recall \cite{ref:ChaNagQueRao-16} that for a measurable map \([\tinit, \tfin]\ni t\mapsto z(t)\in\R^\nu\), the \(\Lp 0\)-norm \(\norm{z}_{\Lp 0([\tinit, \tfin])}\) of \(z\) is defined to be
		\begin{equation}
		\label{e:L0 def}
			\norm{z}_{\Lp 0([\tinit, \tfin])} \Let \Leb\bigl(\{t\in[\tinit, \tfin]\mid z(t) \neq 0\}\bigr);
		\end{equation}
		in other words, \(\norm{z}_{\Lp 0([\tinit, \tfin])}\) is the Lebesgue measure of the support of the map \(z\). It follows, therefore, that
		\begin{equation}
		\label{e:L0 integral}
			\norm{z}_{\Lp 0([\tinit, \tfin])} = \tfin - \int_{\tinit}^{\tfin} \indic{\{0\}}(z(t))\,\dd t.
		\end{equation}

%%%%%%%%%%%%%%%%%%%%%%%%%%%%%%%%%%%%%%%%%%%%%%%%%%%%%%%%%%%%%%%%%%%%%%%%%%%%%%%%
	\section{Multiplexed sparsest reachability}
	\label{s:reach}
%%%%%%%%%%%%%%%%%%%%%%%%%%%%%%%%%%%%%%%%%%%%%%%%%%%%%%%%%%%%%%%%%%%%%%%%%%%%%%%%
		The standard reachability problem in control theory consists of finding, if possible, an admissible control such that, starting from a given initial state, the system evolves under that control to attain a given final state at the end of a given time duration. In the setting of the ensemble \eqref{e:ensemble}, this standard reachability problem admits a natural extension: that of simultaneously transferring the states of each individual member of the ensemble \eqref{e:ensemble} of control systems from given initial to given final states over the time interval \([\tinit, \tfin]\). Moreover, a single control channel must be shared between \(N\) controls, and furthermore, we demand that the individual controls are as sparse as possible.

		Formally, our objective in this section is to \emph{characterize}, for each \(k = 1, \ldots, N\), a control \(\subcon{k}:[\tinit, \tfin]\lra\subadmact{k}\) such that,
		\begin{enumerate}[label=(R-\roman*), leftmargin=*, align=left, widest=iii]
			\item \label{reach:boundary} given initial states \(\substinit{k}\) and final states \(\substfin{k}\) of the \(k\)-th control system in \eqref{e:ensemble}, ensure that \(\subst{k}(\tinit) = \substinit{k}\) and \(\subst{k}(\tfin) = \substfin{k}\),
			\item \label{reach:multiplexed} at a.e.\ \(t\in[\tinit, \tfin]\), at most one \(\subcon{k}(t)\) may be non-zero, and
			\item \label{reach:sparsity} the controls \(\subcon{k}:[\tinit, \tfin]\lra\R^{m_k}\) are set to \(0\) for the longest possible duration.
		\end{enumerate}
		As we shall see momentarily, our characterization of such controls (in Theorem \ref{t:reach}) leads to a certain two-point boundary value problem that can be solved numerically using off-the-shelf numerical routines; to wit, our characterization of the control leads to its computation.

		Given just the requirement \ref{reach:boundary}, if the individual linear control systems in \eqref{e:ensemble} are controllable and the admissible action sets are unconstrained, we know that there exists a control satisfying \ref{reach:boundary}. In such a scenario, it would be natural to define an optimal control problem that minimizes, for instance, an objective function consisting of a convex quadratic function of the individual controls with \ref{reach:boundary} as the set of constraints. Such an unconstrained control problem is classical, and has a well-known solution, that can be obtained, for instance, with the aid of the classical Pontryagin maximum principle \cite[Theorem 22.2]{ref:Cla-13}. However, the resulting control may not satisfy \ref{reach:multiplexed}-\ref{reach:sparsity}. The reachability requirement \ref{reach:boundary} is meaningful if the admissible action sets \(\subadmact{k}\) are not equal to \(\R^{m_k}\); otherwise each individual reachability manoeuvre can be executed in arbitrarily small time with the help of very large control actions provided that the individual systems are controllable, and there is nothing further to do. It is a natural and standard assumption in reachability problems that the admissible action sets \(\subadmact{k}\) are compact for each \(k = 1, \ldots, N\), and we adhere to this assumption throughout this section.

		Since the admissible action sets \((\subadmact{k})_{k=1}^N\) are compact by assumption and \(\tfin\) is given, the stipulations \ref{reach:boundary}-\ref{reach:multiplexed} posed above may be too tight for a feasible control to exist unless \(\tfin\) is sufficiently large. The question of existence of controls satisfing \ref{reach:boundary}-\ref{reach:multiplexed} is tackled by standard minimum time optimal control problems for linear control systems that are very well-understood, and therefore we do not address this question here.%Our result (Theorem \ref{t:reach} below,) is limited to a \emph{characterization} of a control that satisfies \ref{reach:boundary}-\ref{reach:sparsity}.

		Observe that due to the requirement \ref{reach:multiplexed}, whenever a particular control is set to \(0\), the corresponding control system states evolve in open-loop under its unforced/natural dynamics. Consequently, on the one hand, the condition \ref{reach:boundary} cannot be satisfied, in general, by simply concatenating over time \(N\) controls that execute the individual reachability manoeuvres. On the other hand, it may be possible to accomplish the joint reachability manoeuvre in lesser time than the sum of the minimumum times needed for the individual manoeuvres.

		Let us turn to constructing an optimal control problem that takes into account the complete set of requirements \ref{reach:boundary}-\ref{reach:sparsity}. It follows from \ref{reach:boundary}-\ref{reach:sparsity} that the joint system \eqref{e:joint system} must be considered in the process of searching for a set of feasible controls; no decentralized strategy will work. The boundary conditions \ref{reach:boundary} for the individual control systems in \eqref{e:ensemble} lift in an obvious way to boundary conditions for the joint system \eqref{e:joint system}. The device that permits us to encode \ref{reach:multiplexed} for the ensemble \eqref{e:ensemble} into the context of \eqref{e:joint system} is our definition of the set \(\admactreach\) of admissible actions for \eqref{e:joint system} for our reachability problem, given by
		\begin{equation}
		\label{e:admact reach}
			\admactreach \Let \bigcup_{k=1}^N \overset{\text{\(N\)-fold}}{\overbrace{\{0\}\times\cdots\times\subadmact{k}\times\cdots\times\{0\}}} \subset \R^m,
		\end{equation}
		where the set \(\subadmact{k}\) appears as the \(k\)-th factor in the product. Indeed, if for some \(t\in[\tinit, \tfin]\) we have \(\con(t)\in\admactreach\) with \(\subcon{k}(t)\neq 0\), then \(\subcon{\ell}(t) = 0\) for all \(\ell\neq k\). As can be readily checked, the converse also holds. The sparsity requirement \ref{reach:sparsity}, however, is more difficult to ensure, and we do it in an unconventional fashion. We follow the footsteps of \cite{ref:ChaNagQueRao-16}, and consider \(\Lp 0\)-norms of the individual controls in \eqref{e:ensemble} to ensure sparsity, but with an important difference with respect to \cite{ref:ChaNagQueRao-16} that is highlighted in Remark \ref{r:difference}.

		We pick a family \((\lambda_k)_{k=1}^N\) of positive weights as a mechanism to prioritize the set of controls \((\subcon{k})_{k=1}^N\). For each \(k = 1, \ldots, N\), we define \(\lambda_k' \Let \sum_{\ell\in\{1, \ldots, N\}\setmin\{k\}} \lambda_\ell\), and let \(\wt\lambda \Let \sum_{k=1}^n \lambda_k\). In the light of the preceding discussion, we define the optimal control problem:
		\begin{equation}
		\label{e:reach problem prelim}
		\begin{aligned}
			\minimize_{(\subcon{k})_{k=1}^N}	& && \sum_{k=1}^N \lambda_k \norm{\subcon{k}}_{\Lp 0([\tinit, \tfin])}\\
			\sbjto		& && \begin{cases}
				\tdsubst{k}(t) = A_k \subst{k}(t) + B_k \subcon{k}(t),\\
				\subst{k}(\tinit) = \substinit{k}, \quad \subst{k}(\tfin) = \substfin{k},\\
				\text{at most one \(\subcon{k}\) is non-zero at each \(t\in[\tinit, \tfin]\)},\\
				[0, \tfin]\ni t\mapsto \subcon{k}(t)\in\subadmact{k} \text{ measurable},\\
				k = 1, \ldots, N.
			\end{cases}
		\end{aligned}
		\end{equation}
		The second and the third constraints in \eqref{e:reach problem prelim} ensure that \ref{reach:boundary} and \ref{reach:multiplexed} hold, respectively, and the objective function in \eqref{e:reach problem prelim} penalizes the time on which the controls are non-zero, leading to sparsity of the controls complying with \ref{reach:sparsity}. Define the stacked vectors \(\stinit \Let \pmat{\substinit{1}\\\smash[t]{\vdots}\\ \substinit{N}}\) and \(\stfin \Let \pmat{\substfin{1}\\\smash[t]{\vdots}\\\substfin{N}}\). Now \eqref{e:L0 def} and \eqref{e:L0 integral} permit us to rewrite \eqref{e:reach problem prelim} in a more conventional integral form with respect to the joint system \eqref{e:joint system} and the admissible action set defined in \eqref{e:admact reach}:
		\begin{equation}
		\label{e:reach problem}
		\begin{aligned}
			\minimize_{\con}	& && - \int_{\tinit}^{\tfin} \sum_{k=1}^N \lambda_k \indic{\{0\}}(\subcon{k}(t))\,\dd t\\
			\sbjto				& && \begin{cases}
				\text{dynamics \eqref{e:joint system} under the definition \eqref{e:joint state-actions}},\\
				\st(\tinit) = \stinit, \quad \st(\tfin) = \stfin,\\
				[0, \tfin]\ni t\mapsto \con(t)\in\admactreach \text{ measurable}.
			\end{cases}
		\end{aligned}
		\end{equation}
		Clearly, \eqref{e:reach problem} is well-defined due to boundedness of the instantaneous cost function. If a feasible control \(\con\) for \eqref{e:joint system} is supplied by \eqref{e:reach problem}, then projections to appropriate factors yield the individual controls \((\subcon{k})_{k=1}^N\) satisfying the constraints of \eqref{e:reach problem prelim}. If \(\con\opt\) is a minimizer of \eqref{e:reach problem}, with some abuse of terminology, we also say that the state-action trajectory \(t\mapsto \bigl(\st\opt(t), \con\opt(t)\bigr)\) is a \emph{minimizer} of \eqref{e:reach problem}, where \(\st\opt\) is the solution of \eqref{e:joint system} under \(\con\opt\).

		Note that the objective function in \eqref{e:reach problem} contains terms that are discontinuous in the individual controls. This leads to difficulties in solving \eqref{e:reach problem} using standard tools. We employ a nonsmooth Pontryagin maximum principle from \cite[Chapter 22]{ref:Cla-13} to characterize solutions of \eqref{e:reach problem}. Our first main result is the following theorem:

		\begin{theorem}
		\label{t:reach}
			Consider the optimal control problem \eqref{e:reach problem}, and refer to the notations introduced in this section. If \([\tinit, \tfin]\ni t\mapsto \bigl(\st\opt(t), \con\opt(t)\bigr)\in\R^d\times\admactreach\) is a local minimizer of \eqref{e:reach problem}, then there exist a scalar \(\eta = 0\) or \(1\), an absolutely continuous map
			\[
				[\tinit, \tfin]\ni t\mapsto \adj(t) \Let \pmat{\subadj{1}(t)\\\smash[t]{\vdots}\\\subadj{N}(t)}\in\R^{d_1}\times\cdots\times\R^{d_N}
			\]
			and a measurable map
			\[
				[\tinit, \tfin]\ni t\mapsto \reachmux(t)\in\{1, \ldots, N\},
			\]
			such that for a.e.\ \(t\in[\tinit, \tfin]\), either\\
			\(\circ\) \(\eta = 1\) and
			\[
			\left\{
			\begin{aligned}
				& \dot\st\opt(t) = A \st\opt(t) + B \con\opt(t),\quad \st\opt(\tinit) = \stinit, \quad \st\opt(\tfin) = \stfin,\\
				& \dot\adj(t) = -A\trnsp \adj(t),\\
				& \subcon{k}\opt(t) \in \begin{cases}
					\begin{dcases}
					\argmax_{v\in\subadmact{k}}\inprod{B_k\trnsp\subadj{k}(t)}{v}\\
						\qquad \text{if }\max_{v\in\subadmact{k}} \inprod{B_k\trnsp\subadj{k}(t)}{v} > \lambda_k,\\
					\{0\}\cup\argmax_{v\in\subadmact{k}}\inprod{B_k\trnsp\subadj{k}(t)}{v}\\
						\qquad \text{if }\max_{v\in\subadmact{k}} \inprod{B_k\trnsp\subadj{k}(t)}{v} = \lambda_k,\\
					\{0\}\\
						\qquad \text{otherwise}.
				\end{dcases}
					\\
						\qquad\qquad \text{if } k = \reachmux(t),\\
					\{0\}
					\\
						\qquad\qquad \text{otherwise}.
				\end{cases}
				\\
				& \qquad \text{for each } k = 1, \ldots, N,\\
				& \reachmux(t)\in \argmax_{k\in\{1, \ldots, N\}} \max_{v\in\subadmact{k}} \Bigl\{ \inprod{B_k\trnsp \subadj{k}(t)}{v} + \lambda_k \indic{\{0\}}(v) + \lambda_k'\Bigr\},
			\end{aligned}
			\right.
			\]
			or\\
			\(\circ\) \(\eta = 0\) and
			\[
			\left\{
			\begin{aligned}
				& \dot\st\opt(t) = A \st\opt(t) + B \con\opt(t),\quad \st\opt(\tinit) = \stinit, \quad \st\opt(\tfin) = \stfin,\\
				& \dot\adj(t) = -A\trnsp \adj(t),\\
				& \subcon{k}\opt(t) \in \begin{dcases}
					\argmax_{v\in\subadmact{k}}\inprod{B_k\trnsp \subadj{k}(t)}{v}	& \text{if } k = \reachmux(t),\\
					\{0\}	& \text{otherwise}.
				\end{dcases}
				\\
				& \qquad \text{for each } k = 1, \ldots, N,\\
				& \reachmux(t)\in \argmax_{k\in\{1, \ldots, N\}}\max_{v\in\subadmact{k}}\inprod{B_k\trnsp \subadj{k}(t)}{v}.
			\end{aligned}
			\right.
			\]
			Moreover, irrespective of whether \(\eta = 1\) or \(0\), the map
			\[
				[0, \tfin]\ni t\mapsto \inprod{\adj(t)}{A\st\opt(t) + B\con\opt(t)} + \eta \sum_{k=1}^N \lambda_k \indic{\{0\}}(\subcon{k}\opt(t))\in\R
			\]
			is a constant a.e.
		\end{theorem}

		\begin{remark}
		\label{r:reach description}
			\textup{Theorem \ref{t:reach} gives a characterization of the optimal control \(\con\opt\) in the same spirit as Euler's first order necessary condition for a minimum: if \(O\subset\R^\nu\) is a non-empty open set, \(g:O\lra\R\) is a continuously differentiable function, and \(z\opt\in O\) is a minimizer of \(g\), then \(\nabla g(z\opt) = 0\). Armed with such a characterization, one lets algorithms find the optimizers. In the case of Theorem \ref{t:reach}, observe that the characterization in both the cases of \(\eta = 1\) and \(\eta = 0\) leads to a boundary value problem consisting of a family of \(2d\) scalar differential equations with \(2d\) boundary values; as such they are well-posed problems. One typically employs shooting methods to solve such boundary value problems, and there are many such techniques available today, see, e.g., \cite{ref:Ros-15} for details. The multiplexer that we want is given by the map \(\reachmux\).}
		\end{remark}

		\begin{remark}
			\textup{As with any optimal control problem, the question of existence of a minimizer in \eqref{e:reach problem} arises naturally. Unfortunately, such an existential result is not known at present. Classical existential results given in, e.g., \cite{ref:FleRis-75}, do not apply to \eqref{e:reach problem} due to discontinuities in the instantaneous cost function there; for the same reason, classical techniques for proving existence of minimizers also cease to apply in the context of \eqref{e:reach problem}.}
		\end{remark}

		\begin{remark}
		\label{r:difference}
			\textup{A direct translation of the \(\Lp 0\)-cost considered in \cite{ref:ChaNagQueRao-16} to our context would be the \(\Lp 0\)-cost \(\norm{\con}_{\Lp 0([\tinit, \tfin])}\) of the joint control \(\con \Let \pmat{\subcon{1}\\\smash[t]{\vdots}\\\subcon{N}}\). Here, instead, we work with a (weighted) sum \(\sum_{k=1}^N \lambda_k \norm{\subcon{k}}_{\Lp 0([\tinit, \tfin])}\)  of \(\Lp 0\)-costs of the individual controls \(\subcon{k}\). Both the costs minimize the duration of time that the controls are non-zero. In the former case, the instantaneous cost is positive only if \emph{all} the individual controls are precisely \(0\). In the later case, the instaneous cost increases on every subset of \([\tinit, \tfin]\) of positive measure whenever \emph{at least one} of the individual controls attains the value \(0\). Cf.\ \cite[Remark 1]{ref:ChaNagQueRao-16}.}
		\end{remark}

		The following Corollary isolates the important case of each \(\subcon{k}\) being \([-1, 1]\)-valued:

		\begin{corollary}
		\label{c:reach scalar}
			Consider the optimal control problem \eqref{e:reach problem} with each \(\subadmact{k} = [-1, 1]\), and refer to the notations introduced in this section. If \([\tinit, \tfin]\ni t\mapsto \bigl(\st\opt(t), \con\opt(t)\bigr)\in\R^d\times\admactreach\) is a local minimizer of \eqref{e:reach problem}, then then there exist a scalar \(\eta = 0\) or \(1\), an absolutely continuous map
			\[
				[\tinit, \tfin]\ni t\mapsto \adj(t) \Let \pmat{\subadj{1}(t)\\\smash[t]{\vdots}\\\subadj{N}(t)}\in\R^{d_1}\times\cdots\times\R^{d_N}
			\]
			and a measurable map
			\[
				[\tinit, \tfin]\ni t\mapsto \reachmux(t)\in\{1, \ldots, N\},
			\]
			such that for a.e.\ \(t\in[\tinit, \tfin]\), either\\
			\(\circ\) \(\eta = 1\) and
			\[
			\left\{
			\begin{aligned}
				& \dot\st\opt(t) = A \st\opt(t) + B \con\opt(t),\quad \st\opt(\tinit) = \stinit, \quad \st\opt(\tfin) = \stfin,\\
				& \dot\adj(t) = -A\trnsp \adj(t),\\
				& \subcon{k}\opt(t) \in \begin{cases}
					\begin{dcases}
					\sgn\bigl(B_k\trnsp\subadj{k}(t)\bigr)				& \text{if } \abs{B_k\trnsp\subadj{k}(t)} > \lambda_k,\\
					\{0\}\cup\sgn\bigl(B_k\trnsp\subadj{k}(t)\bigr)	& \text{if } \abs{B_k\trnsp\subadj{k}(t)} = \lambda_k,\\
					\{0\}	& \text{otherwise}.
				\end{dcases}
					\\
						\qquad\qquad \text{if } k = \reachmux(t),\\
					\{0\}
					\\
						\qquad\qquad \text{otherwise}.
				\end{cases}
				\\
				& \qquad \text{for each } k = 1, \ldots, N,\\
				& \reachmux(t)\in \argmax_{k\in\{1, \ldots, N\}} \begin{cases}
					        \lambda_k'  +  \abs{B_k\trnsp \subadj{k}(t)}        & \text{if }\abs{B_k\trnsp \subadj{k}(t)}    \geq   \lambda_k,\\    %%%% to be edited
				               \wt\lambda  	&        \text{otherwise}.
				\end{cases}
					%\max_{v\in[-1, 1]} \Bigl\{ B_k\trnsp \subadj{k}(t)v + \lambda_k \indic{\{0\}}(v) + \lambda_k'\Bigr\},
			\end{aligned}
			\right.
			\]
			or\\
			\(\circ\) \(\eta = 0\) and
			\[
			\left\{
			\begin{aligned}
				& \dot\st\opt(t) = A \st\opt(t) + B \con\opt(t),\quad \st\opt(\tinit) = \stinit, \quad \st\opt(\tfin) = \stfin,\\
				& \dot\adj(t) = -A\trnsp \adj(t),\\
				& \subcon{k}\opt(t) \in \begin{dcases}
					\sgn\bigl(B_k\trnsp \subadj{k}(t)\bigr)	& \text{if } k = \reachmux(t),\\
					\{0\}	& \text{otherwise},
				\end{dcases}
				\\
				& \qquad \text{for each } k = 1, \ldots, N,\\
				& \reachmux(t)\in \argmax_{k\in\{1, \ldots, N\}}\abs{B_k\trnsp \subadj{k}(t)}.
			\end{aligned}
			\right.
			\]
		\end{corollary}

%%%%%%%%%%%%%%%%%%%%%%%%%%%%%%%%%%%%%%%%%%%%%%%%%%%%%%%%%%%%%%%%%%%%%%%%%%%%%%%%
	\section{Multiplexed sparse LQ control}
	\label{s:LQ}
%%%%%%%%%%%%%%%%%%%%%%%%%%%%%%%%%%%%%%%%%%%%%%%%%%%%%%%%%%%%%%%%%%%%%%%%%%%%%%%%
		In this section we address the question of controlling the ensemble of linear control systems \eqref{e:ensemble} by minimizing a standard quadratic instantaneous cost on the individual states and controls while enforcing the constraint that at any instant of time only one of the controls is ``active,'' and demanding that the individual controls are sparse.

		Formally, our objective is to \emph{characterize}, for each \(k = 1, \ldots, N\), a control \(\subcon{k}:[\tinit, \tfin]\lra\R^{m_k}\) such that, given initial states \(\substinit{k}\) of the \(k\)-th control system in \eqref{e:ensemble},
		\begin{enumerate}[label=(LQ-\roman*), leftmargin=*, align=left, widest=iii]
			\item \label{LQ:cost} \(\subcon{k}\) minimizes a standard quadratic objective function with respect to the states and the controls of the \(k\)-th system in \eqref{e:ensemble},
			\item \label{LQ:multiplexed} at a.e.\ \(t\in[\tinit, \tfin]\), at most one \(\subcon{k}(t)\) may be non-zero, and
			\item \label{LQ:sparsity} the controls \(\subcon{k}:[\tinit, \tfin]\lra\R^{m_k}\) are set to \(0\) whenever possible.
		\end{enumerate}

		The standard linear quadratic control problem for individual members of the ensemble \eqref{e:ensemble} consists of the following: For each \(k = 1, \ldots, N\), let symmetric and non-negative definite matrices \(Q_k, \hat Q_k\in\R^{d_k\times d_k}\) and a symmetric and positive definite matrix \(R_k\in\R^{m_k\times m_k}\) be given, and let \(\substinit{k}\in\R^{d_k}\) be a given initial state. Consider the following optimal control problem for each \(k = 1, \ldots, N\):
		\begin{equation}
		\label{e:LQ problem individual}
		\begin{aligned}
			\minimize_{\subcon{k}}	&&& \tfrac{1}{2} \int_{\tinit}^{\tfin} \left( \inprod{\subst{k}(t)}{Q_k \subst{k}(t)} + \inprod{\subcon{k}(t)}{R_k \subcon{k}(t)} \right)\,\dd t\\
									&&& + \tfrac{1}{2} \inprod{\subst{k}(\tfin)}{\hat Q_k \subst{k}(\tfin)}\\
			\sbjto					&&& \begin{cases}
				\tdsubst{k}(t) = A_k \subst{k}(t) + B \subcon{k}(t),\\
				\subst{k}(\tinit) = \substinit{k},\\
				[\tinit, \tfin]\ni t\mapsto \subcon{k}(t)\in\R^{m_k}\text{ measurable}.
			\end{cases}
		\end{aligned}
		\end{equation}
		The classical theory of linear quadratic optimal control \cite[Chapter 6]{ref:Lib-12} asserts that a solution of \eqref{e:LQ problem individual} exists under the preceding conditions, determined completely by the so-called Riccati differential equation \cite[Equation (6.14)]{ref:Lib-12}
		\[
		\left\{
		\begin{aligned}
			& \dot P_k(t) = - A_k\trnsp P_k(t) - P_k(t) A_k - Q_k + P_k(t) B_k R_k\inverse B_k\trnsp P_k(t),\\
			& \qquad\qquad\text{for a.e.\ }t\in[\tinit, \tfin],\\
			& P_k(\tfin) = \hat Q_k,
		\end{aligned}
		\right.
		\]
		with the (unique) optimal state-action trajectory \([\tinit, \tfin]\ni t\mapsto \bigl(\subst{k}\opt(t), \subcon{k}\opt(t)\bigr)\in\R^{d_k}\times\R^{m_k}\) expressed in terms of \(P_k\) as
		\[
		\left\{
		\begin{aligned}
			& \tdsubst{k}\opt(t) = A_k \subst{k}\opt(t) + B_k \subcon{k}\opt(t) \quad \text{for a.e.\ }t\in[\tinit, \tfin],\\
			& \subcon{k}\opt(t) = R_k\inverse B_k\trnsp P_k(t) \subst{k}\opt(t) \quad \text{for a.e.\ }t\in[\tinit, \tfin],\\
			& \subst{k}(\tinit) = \substinit{k}.
		\end{aligned}
		\right.
		\]
		Of course, the conditions \ref{LQ:multiplexed}-\ref{LQ:sparsity} cannot be ensured by simply solving the individual LQ problems \eqref{e:LQ problem individual} because neither the multiplexing requirement \ref{LQ:multiplexed} nor the sparsity desideratum \ref{LQ:sparsity} is incorporated into the problem \eqref{e:LQ problem individual} since it is defined separately for each \(k\).

		Let us construct an optimal control problem that accounts for the complete set of requirements \ref{LQ:cost}-\ref{LQ:sparsity}. It follows from \ref{LQ:multiplexed}-\ref{LQ:sparsity} that the joint system \eqref{e:joint system} must be considered in the process of searching for a feasible set of controls. The device that permits us to encode \ref{LQ:multiplexed} into the context of \eqref{e:joint system} is the definition of the set \(\admactLQ\) of admissible actions for \eqref{e:joint system} given by
		\begin{equation}
		\label{e:admact LQ}
			\admactLQ \Let \bigcup_{k=1}^N \overset{\text{\(N\)-fold}}{\overbrace{\{0\}\times\cdots\times\R^{m_k}\times\cdots\times\{0\}}} \subset \R^m,
		\end{equation}
		where \(\R^{m_k}\) appears as the \(k\)-th factor in the product. (The definition of \(\admactLQ\) mirrors that of \(\admactreach\) defined in \eqref{e:admact reach}.) The set \(\admactLQ\) is a non-convex cone and star-shaped about \(0\in\R^m\).\footnote{Recall that a set \(C\subset\R^\nu\) is a \emph{cone} if for every \(z\in C\) and \(\alpha \ge 0\), the point \(\alpha z\) belongs to \(C\), and \(C\) is \emph{star shaped about \(0\)} if for every \(y\in C\), the straight line segment joining \(0\) to \(y\) is contained in \(C\).} Observe that a well-defined performance index is already present in \ref{LQ:cost}; it is therefore not possible to simultaneously stipulate \emph{maximal} sparsity in the controls, for that would lead to two different performance indices.\footnote{The issue of pareto optimality, while interesting in our context, is not treated in this article.} Instead we enforce sparsity by \(\Lp 0\)-regularizing the individual cost functions relative to the corresponding controls along the lines of \cite{ref:SriCha-16}; the \(\Lp 0\)-regularization parameters influence the extent of sparsity that arise as a consequence.

		More formally, let \(\lambda_k > 0\) be a regularization parameter for each \(k = 1, \ldots, N\), let \(\lambda_k' \Let \sum_{\ell\in\{1, \ldots, N\}\setmin\{k\}} \lambda_\ell\) for each \(k\), and let \(\wt\lambda \Let \sum_{k=1}^N \lambda_k\). Define the stacked vector \(\stinit \Let \pmat{\substinit{1}\\\smash[t]{\vdots}\\ \substinit{N}}\), block-diagonal matrices \(Q \Let \blkdiag(Q_1, \ldots, Q_N)\), \(R \Let \blkdiag(R_1, \ldots, R_N)\), \(\hat Q \Let \blkdiag(\hat Q_1, \ldots, \hat Q_N)\). In view of the preceding discussion, \eqref{e:L0 def}, and \eqref{e:L0 integral}, we arrive at the following optimal control problem in an integral form:
		\begin{equation}
		\label{e:LQ problem}
		\begin{aligned}
			\minimize_{\con}	&&& \int_{\tinit}^{\tfin} \Biggl(\tfrac{1}{2}\bigl(\inprod{\st(t)}{Q \st(t)} + \inprod{\con(t)}{R \con(t)}\bigr)\\
								&&& \qquad - \sum_{k=1}^N \lambda_k \indic{\{0\}}(\subcon{k}(t)) \Biggr)\,\dd t + \tfrac{1}{2}\inprod{\st(\tfin)}{\hat Q \st(\tfin)}\\
			\sbjto				&&& \begin{cases}
				\text{dynamics \eqref{e:joint system} under the definition \eqref{e:joint state-actions}},\\
				\st(\tinit) = \stinit,\\
				[\tinit, \tfin]\ni t\mapsto u(t)\in\admactLQ\text{ measurable}.
			\end{cases}
		\end{aligned}
		\end{equation}
		Quite clearly, \eqref{e:LQ problem} is well-defined due to positive definiteness of the matrix \(R\). If a feasible control \(\con\) for \eqref{e:joint system} is supplied by \eqref{e:LQ problem}, then projections to appropriate factors yield the individual controls \((\subcon{k})_{k=1}^N\) satisfying the constraints of \eqref{e:LQ problem individual}. As in the case of \eqref{e:reach problem}, if \(\con\opt\) is a minimizer of \eqref{e:LQ problem}, with some abuse of terminology, we also say that the state-action trajectory \(t\mapsto \bigl(\st\opt(t), \con\opt(t)\bigr)\) is a \emph{minimizer} of \eqref{e:reach problem}, where \(\st\opt\) is the solution of \eqref{e:joint system} under \(\con\opt\).

		Note that the objective function in \eqref{e:LQ problem} contains terms that are discontinuous in the individual controls, which leads to difficulties in solving \eqref{e:LQ problem} using standard tools. We employ a nonsmooth Pontryagin maximum principle from \cite[Chapter 22]{ref:Cla-13} to characterize solutions of \eqref{e:LQ problem}, and this characterization is the subject of the following theorem:

		\begin{theorem}
		\label{t:LQ}
			Consider the optimal control problem \eqref{e:LQ problem} and refer to the notations introduced in this section. If \([\tinit, \tfin]\ni t\mapsto\bigl(\st\opt(t), \con\opt(t)\bigr)\in\R^d\times\admactLQ\) is a local minimizer of \eqref{e:LQ problem}, then there exist an absolutely continuous map
			\[
				[\tinit, \tfin]\ni t\mapsto \adj(t) \Let \pmat{\subadj{1}(t)\\\smash[t]{\vdots}\\\subadj{N}(t)}\in\R^{d_1}\times\cdots\times\R^{d_N}
			\]
			and a measurable map
			\[
				[\tinit, \tfin]\ni t\mapsto \LQmux(t)\in\{1, \ldots, N\},
			\]
			such that for a.e.\ \(t\in[\tinit, \tfin]\),
			\[
			\left\{
			\begin{aligned}
				& \dot\st\opt(t) = A \st\opt(t) + B \con\opt(t),\quad \st\opt(\tinit) = \stinit,\\
				& \dot\adj(t) = -A\trnsp \adj(t) + Q \st\opt(t),\quad \adj(\tfin) = -\hat Q\st\opt(\tfin),\\
				& \subcon{k}\opt(t) \in \begin{cases}
					\begin{cases}
						\bigl\{R_k\inverse B_k\trnsp \subadj{k}(t)\bigr\}			& \text{if } \norm{\subadj{k}(t)}_{B_k R_k\inverse B_k\trnsp}^2 > 2\lambda_k,\\
						\{0\}\cup\bigl\{R_k\inverse B_k\trnsp \subadj{k}(t)\bigr\}	& \text{if } \norm{\subadj{k}(t)}_{B_k R_k\inverse B_k\trnsp}^2 = 2\lambda_k,\\
						\{0\}	& \text{otherwise},
					\end{cases}
					\\
						\qquad\qquad \text{if } k = \LQmux(t),\\
					\{0\}
					\\
						\qquad\qquad \text{otherwise}.
				\end{cases}
				\\
				& \qquad \text{for each } k = 1, \ldots, N,\\
				& \LQmux(t)\in\\
				& \argmax_{k\in\{1, \ldots, N\}} \begin{cases}
						\lambda_k' + \frac{1}{2}\norm{\subadj{k}(t)}_{B_k R_k\inverse B_k\trnsp}^2	& \text{if }\norm{\subadj{k}(t)}_{B_k R_k\inverse B_k\trnsp}^2 \ge  2\lambda_k,\\
						\wt\lambda	& \text{otherwise}.
					\end{cases}
			\end{aligned}
			\right.
			\]
			Moreover, the map
			\begin{align*}
				[0, \tfin]\ni t\mapsto & \inprod{\adj(t)}{A\st\opt(t) + B\con\opt(t)} + \sum_{k=1}^N \lambda_k \indic{\{0\}}(\subcon{k}\opt(t))\\
				& - \tfrac{1}{2}\bigl(\inprod{\st\opt(t)}{Q \st\opt(t)} + \inprod{\con\opt(t)}{R \con\opt(t)}\bigr)\in\R
			\end{align*}
			is a constant a.e.
		\end{theorem}

		\begin{remark}
		\label{r:LQ description}
			\textup{Theorem \ref{t:LQ} leads to a boundary value problem consisting of a family of \(2d\) scalar differential equations with \(2d\) boundary values; as such it is a well-posed problem. The map \(\LQmux\) is the multiplexer that we want.}
		\end{remark}

		\begin{remark}
			\textup{Just as in \eqref{e:reach problem}, the question of existence of a minimizer in \eqref{e:LQ problem} is natural, and once more, due to discontinuities in the instantaneous cost function in \eqref{e:LQ problem}, classical results and techniques dealing with existence of minimizers do not apply.}
		\end{remark}

%%%%%%%%%%%%%%%%%%%%%%%%%%%%%%%%%%%%%%%%%%%%%%%%%%%%%%%%%%%%%%%%%%%%%%%%%%%%%%%%
	\section{Multiplexed sparse Mayer problem}
	\label{s:Mayer}
%%%%%%%%%%%%%%%%%%%%%%%%%%%%%%%%%%%%%%%%%%%%%%%%%%%%%%%%%%%%%%%%%%%%%%%%%%%%%%%%
		The special case of \(Q = 0\) and \(R = 0\) in \eqref{e:LQ problem} is interesting in its own right; it corresponds to what is commonly known as the multiplexed sparsest Mayer problem. It deserves to be treated separately because the hypotheses of Theorem \ref{t:LQ} do not hold in this setting --- in particular, \(R = 0\) here. 

		Formally, our objective is to \emph{characterize}, for each \(k = 1, \ldots, N\), a control \(\subcon{k}:[\tinit, \tfin]\lra\R^{m_k}\) such that, given initial states \(\substinit{k}\) of the \(k\)-th control system in \eqref{e:ensemble},
		\begin{enumerate}[label=(M-\roman*), leftmargin=*, align=left, widest=iii]
			\item \label{Mayer:cost} \(\subcon{k}\) minimizes a standard quadratic terminal cost function with respect to the states of the \(k\)-th system in \eqref{e:ensemble},
			\item \label{Mayer:multiplexed} at a.e.\ \(t\in[\tinit, \tfin]\), at most one \(\subcon{k}(t)\) may be non-zero, and
			\item \label{Mayer:sparsity} the controls \(\subcon{k}:[\tinit, \tfin]\lra\R^{m_k}\) are set to \(0\) whenever possible.
		\end{enumerate}
		Notice that for the problem to be well-posed, the admissible action sets must be bounded since there is no cost on the control; accordingly we stipulate that the individual admissible action sets \((\subadmact{k})_{k=1}^N\) are all compact with \(0\in\subadmact{k}\) for each \(k = 1, \ldots, N\).

		Let us construct an optimal control problem that accounts for the complete set of requirements \ref{Mayer:cost}-\ref{Mayer:sparsity}. As in the case of \eqref{e:reach problem}, the device that permits us to encode \ref{Mayer:multiplexed} into the context of \eqref{e:joint system} is the definition of the set \(\admactMayer\) of admissible actions for \eqref{e:joint system} given by
		\begin{equation}
		\label{e:admact Mayer}
			\admactMayer \Let \bigcup_{k=1}^N \overset{\text{\(N\)-fold}}{\overbrace{\{0\}\times\cdots\times\subadmact{k}\times\cdots\times\{0\}}} \subset \R^m,
		\end{equation}
		where \(\subadmact{k}\) appears as the \(k\)-th factor in the product. We stipulate maximal sparsity of the controls by placing the \(\Lp 0\)-norms of the individual controls.

		To be precise, for each \(k = 1, \ldots, N\), let \(\stfin_k\in\R^{d_k}\), let \(\hat Q_k\) denote a symmetric and non-negative definite matrix, and let \(\stinit^k\in\R^{d_k}\) denote a given initial state. Let \(\lambda_k > 0\) be a weight for each \(k = 1, \ldots, N\), let \(\lambda_k' \Let \sum_{\ell\in\{1, \ldots, N\}\setmin\{k\}} \lambda_k\), and let \(\wt\lambda \Let \sum_{k=1}^N \lambda_k\). Define the stacked vectors \(\stinit \Let \pmat{\substinit{1}\\\smash[t]{\vdots}\\ \substinit{N}}\), \(\stfin \Let \pmat{\substfin{1}\\\smash[t]{\vdots}\\\substfin{N}}\), and the block-diagonal matrix \(\hat Q \Let \blkdiag\{\hat Q_1, \ldots, \hat Q_N\}\). Consider the problem
		\begin{equation}
		\label{e:Mayer}
		\begin{aligned}
			\minimize_{\con}	&&& - \int_{\tinit}^{\tfin} \sum_{k=1}^N \lambda_k \indic{\{0\}}(\subcon{k}(t))\,\dd t\\
								&&& \qquad + \tfrac{1}{2}\inprod{\st(\tfin)-\stfin}{\hat Q \bigl(\st(\tfin) - \stfin\bigr)}\\
			\sbjto				&&& \begin{cases}
				\text{dynamics \eqref{e:joint system} under the definition \eqref{e:joint state-actions}},\\
				\st(\tinit) = \stinit,\\
				[\tinit, \tfin]\ni t\mapsto u(t)\in\admactMayer\text{ measurable}.
			\end{cases}
		\end{aligned}
		\end{equation}

		Notice that any solution of \eqref{e:Mayer} is maximally sparse while minimizing the terminal cost function; the latter forces small deviations of \(\st\opt(\tfin)\) from the given final state \(\stfin\).

		\begin{remark}
			\textup{Sometimes the problem \eqref{e:Mayer} is employed to achieve approximate reachability by picking a `large' terminal cost. Intuition suggests that if the terminal cost in \eqref{e:Mayer} is large (i.e., the matrix \(\hat Q\) is large), then the resulting control is such that the separation between \(\st\opt(\tfin)\) and \(\stfin\) is small, leading to approximate reachability. The presence of a quadratic term in the terminal cost improves the behaviour of numerical algorithms that are typically employed for multiple shooting methods to solve the two point boundary value problems associated with \eqref{e:Mayer} in comparison to \eqref{e:reach problem}. However, notice that this way sparsity may be compromised to arrive at a smaller terminal cost, and therefore, the two problems \eqref{e:reach problem} and \eqref{e:Mayer} are intrinsically different. Only if the terminal cost in \eqref{e:Mayer} is set to \(0\) if \(\st(\tfin) = \stfin\) and \(+\infty\) otherwise, then the Mayer problem \eqref{e:Mayer} becomes equivalent to the reachability problem \eqref{e:reach problem}.}
		\end{remark}

		The following theorem characterizes solutions of \eqref{e:Mayer}:
		\begin{theorem}
		\label{t:Mayer}
		Consider the optimal control problem \eqref{e:LQ problem}, refer to the notations introduced in this section, and assume that \(Q = 0\) and \(R = 0\). If \([\tinit, \tfin]\ni t\mapsto\bigl(\st\opt(t), \con\opt(t)\bigr)\in\R^d\times\admactMayer\) is a local minimizer of \eqref{e:Mayer}, then there exist an absolutely continuous map
			\[
				[\tinit, \tfin]\ni t\mapsto \adj(t) \Let \pmat{\subadj{1}(t)\\\smash[t]{\vdots}\\\subadj{N}(t)}\in\R^{d_1}\times\cdots\times\R^{d_N}
			\]
			and a measurable map
			\[
				[\tinit, \tfin]\ni t\mapsto \Mayermux(t)\in\{1, \ldots, N\},
			\]
			such that for a.e.\ \(t\in[\tinit, \tfin]\),
			\[
			\left\{
			\begin{aligned}
				& \dot\st\opt(t) = A \st\opt(t) + B \con\opt(t),\quad \st\opt(\tinit) = \stinit,\\
				& \dot\adj(t) = -A\trnsp \adj(t),\quad \adj(\tfin) = -\hat Q\bigl(\st\opt(\tfin) - \stfin\bigr),\\
			& \subcon{k}\opt(t) \in 
			\begin{dcases}
				\begin{dcases}
					\argmax_{v\in\subadmact{k}} \inprod{B_k\trnsp\subadj{k}(t)}{v}\\
					\qquad \text{if } \max_{v\in\subadmact{k}} \inprod{B_k\trnsp\subadj{k}(t)}{v} > \lambda_k,\\
					\{0\}\cup\argmax_{v\in\subadmact{k}} \inprod{B_k\trnsp\subadj{k}(t)}{v}\\
					\qquad \text{if } \max_{v\in\subadmact{k}} \inprod{B_k\trnsp\subadj{k}(t)}{v} = \lambda_k,\\
					\{0\}\\
					\qquad \text{otherwise},
				\end{dcases}
				\\
				\qquad\qquad \text{if } k = \Mayermux(t),\\
				\{0\}
				\\
				\qquad\qquad\text{otherwise},
			\end{dcases}
			\\
			& \qquad\text{for }k = 1, \ldots, N,\\
				& \Mayermux(t)\in \argmax_{k\in\{1, \ldots, N\}} \max_{v\in\subadmact{k}} \Bigl\{\inprod{B_k\trnsp \subadj{k}(t)}{v} + \lambda_k \indic{\{0\}}(v) + \lambda_k'\Bigr\}.
			\end{aligned}
			\right.
			\]
			Moreover, the map
			\begin{align*}
				[0, \tfin]\ni t\mapsto \inprod{\adj(t)}{A\st\opt(t) + B\con\opt(t)} + \sum_{k=1}^N \lambda_k \indic{\{0\}}(\subcon{k}\opt(t))\in\R
			\end{align*}
			is a constant a.e.
		\end{theorem}
		
		\begin{remark}
		\label{r:Mayer description}
			\textup{Theorem \ref{t:Mayer} leads to a boundary value problem consisting of a family of \(2d\) scalar differential equations with \(2d\) boundary values; as such it is a well-posed problem. The map \(\Mayermux\) is our desired multiplexer.}
		\end{remark}

		\begin{remark}
			\textup{Yet again, due to discontinuities in the instantaneous cost function in \eqref{e:Mayer}, classical results and techniques dealing with existence of minimizers do not apply, and the question of existence of minimizers in \eqref{e:Mayer} remains open.}
		\end{remark}

%%%%%%%%%%%%%%%%%%%%%%%%%%%%%%%%%%%%%%%%%%%%%%%%%%%%%%%%%%%%%%%%%%%%%%%%%%%%%%%%
	\section{Numerical experiments}
	\label{s:examples}
\begin{figure}[h]
			\centering
			\includegraphics[scale=0.7]{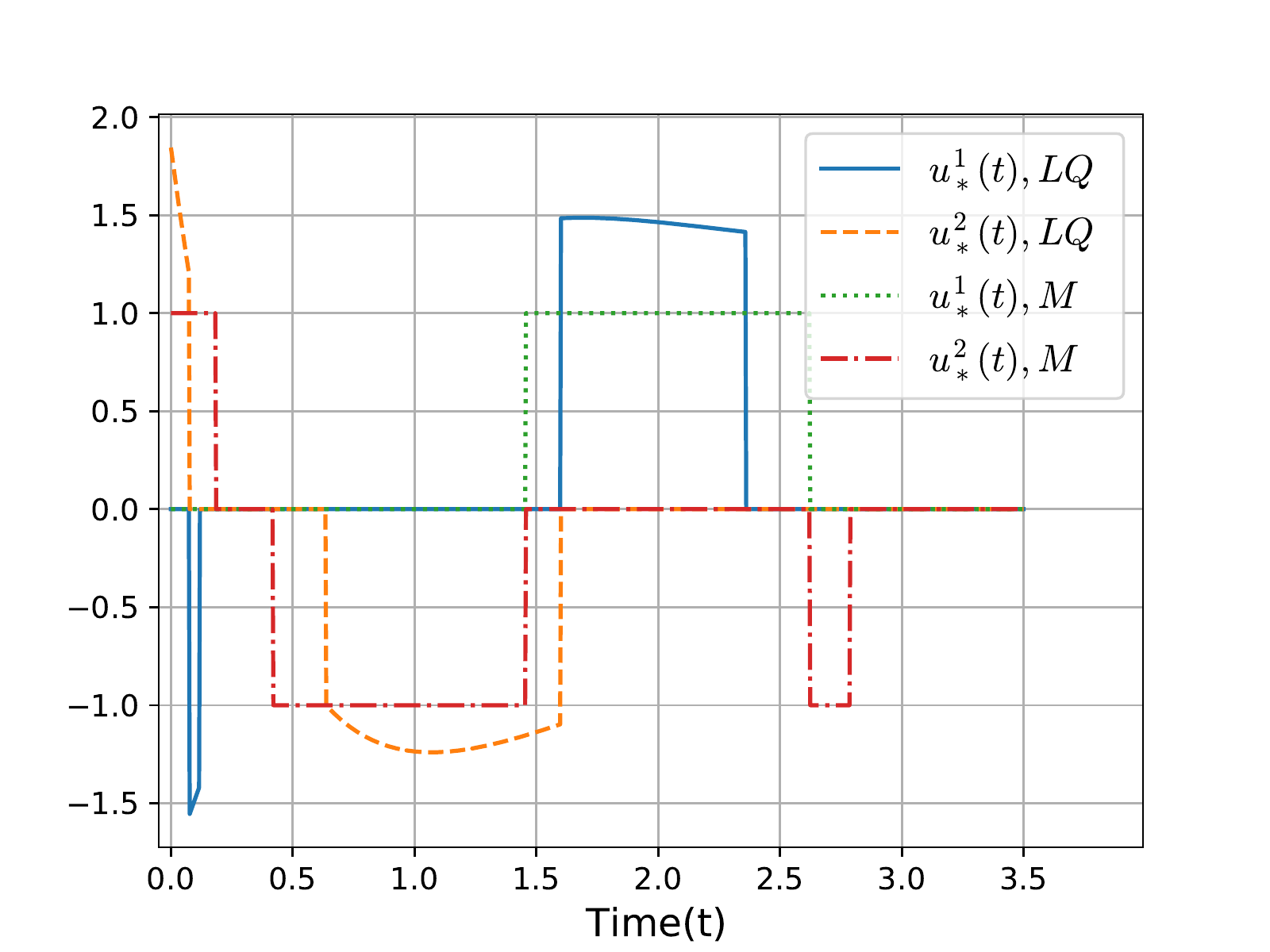}
			\caption{Multiplexed controls for the LQ and the Mayer problems.}
			\label{fig:fig01}
		\end{figure}

\begin{figure}[h]
			\includegraphics[scale=0.7]{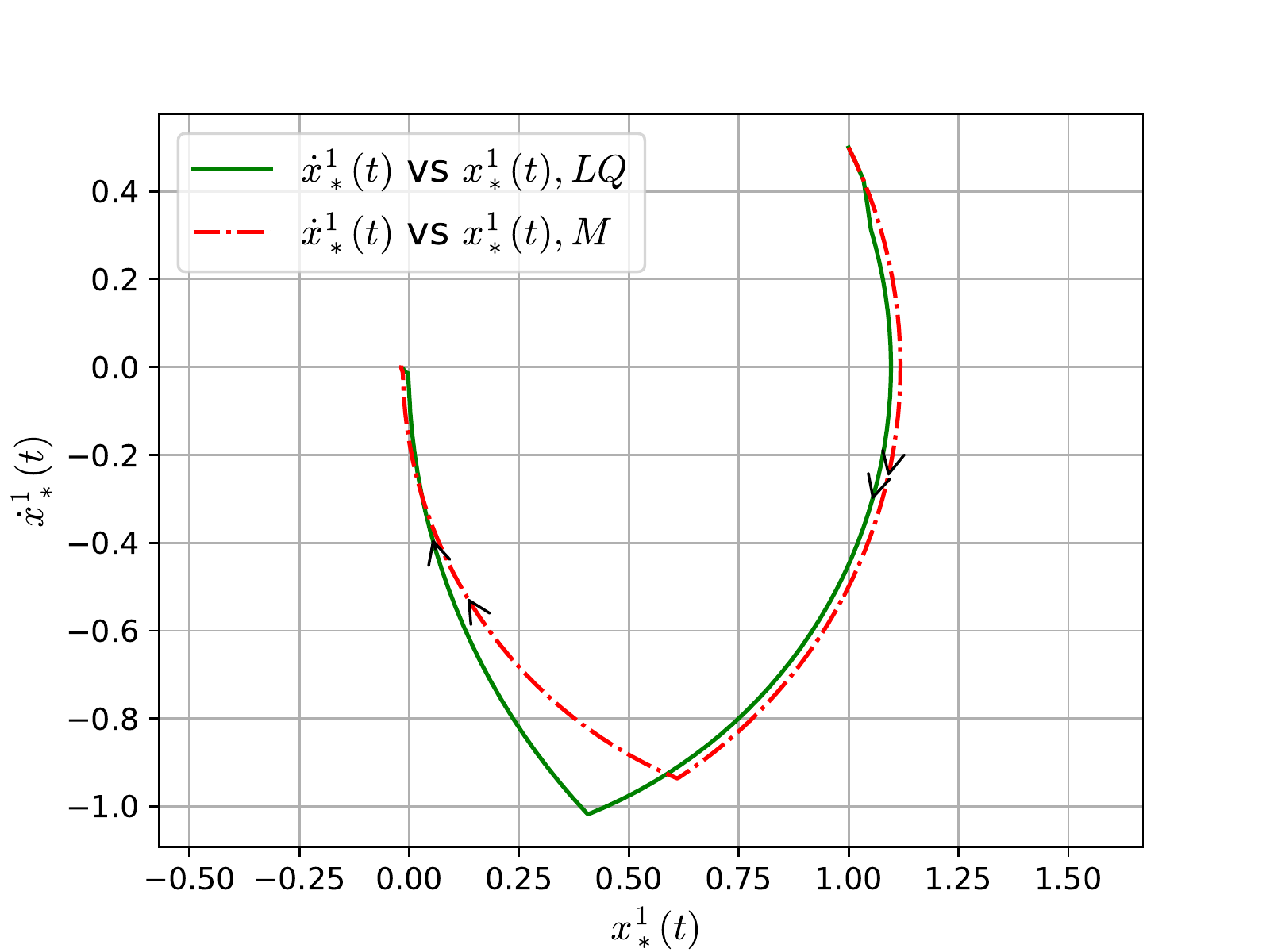}
			\caption{Phase portraits for the LQ and the Mayer problems.}
			\label{fig:fig02}
		\end{figure}

\begin{figure}[h]
			\centering
			\includegraphics[scale=0.7]{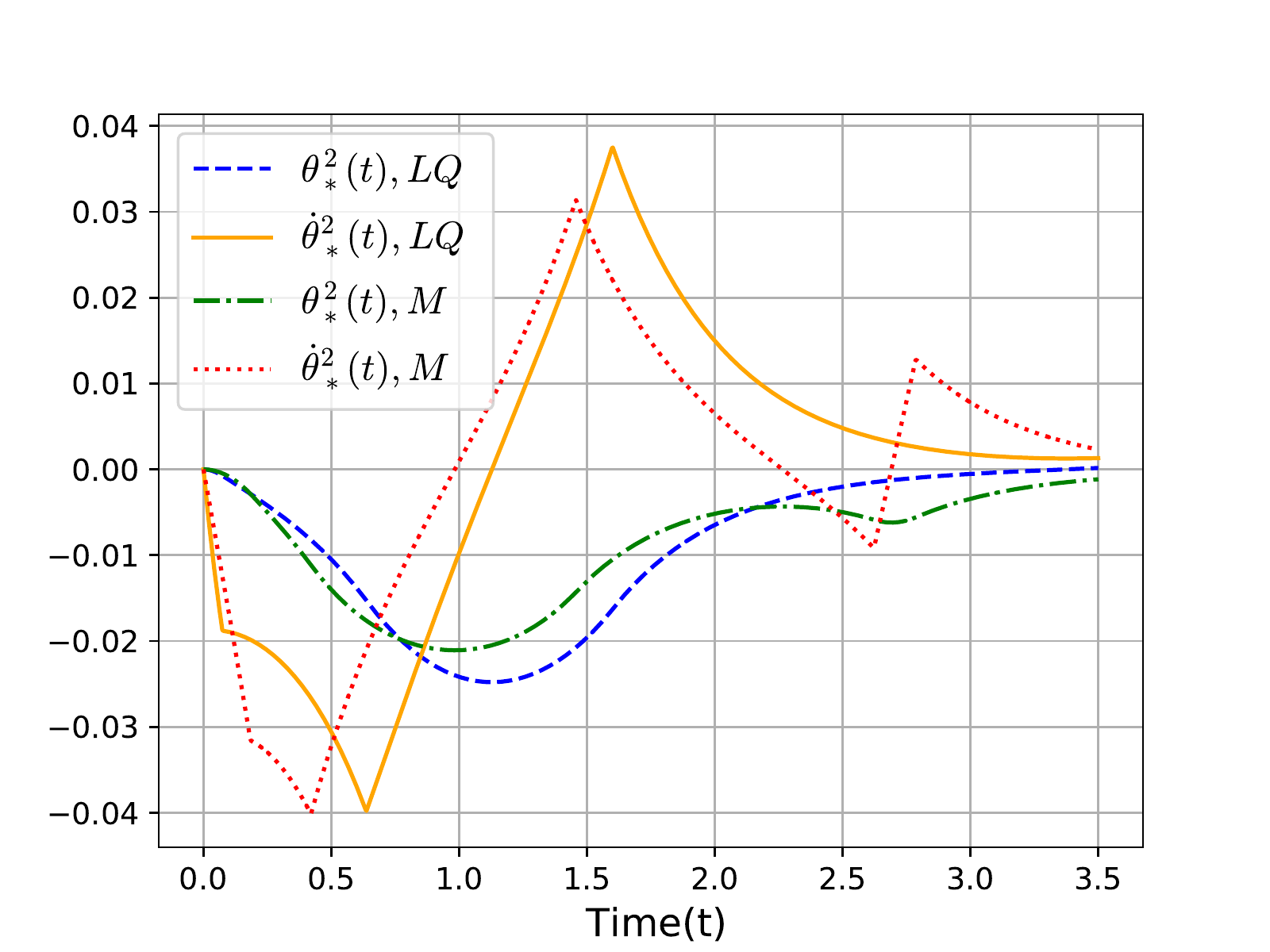}
			\caption{$\theta_*(\cdot) ,\dot{\theta}_*(\cdot)$, for the LQ and the Mayer problems.}
			\label{fig:fig03}
		\end{figure}

\begin{figure}[h]
			\centering
			\includegraphics[scale=0.7]{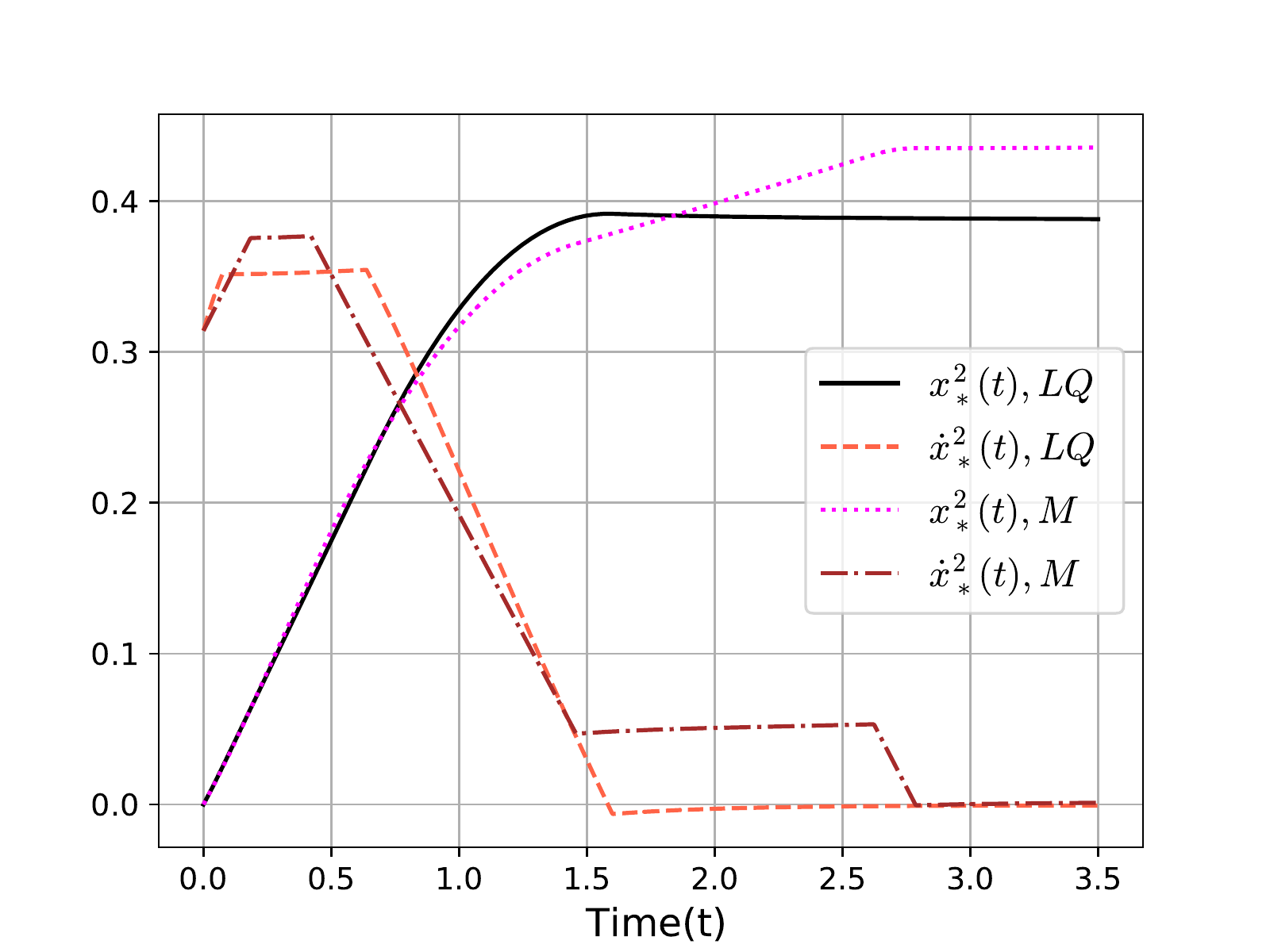}
			\caption{$x_*(\cdot), \dot{x}_*(\cdot)$ for the LQ and the Mayer problems.}
			\label{fig:fig04}
		\end{figure}

\begin{figure}[h]
			\centering
			\includegraphics[scale=0.7]{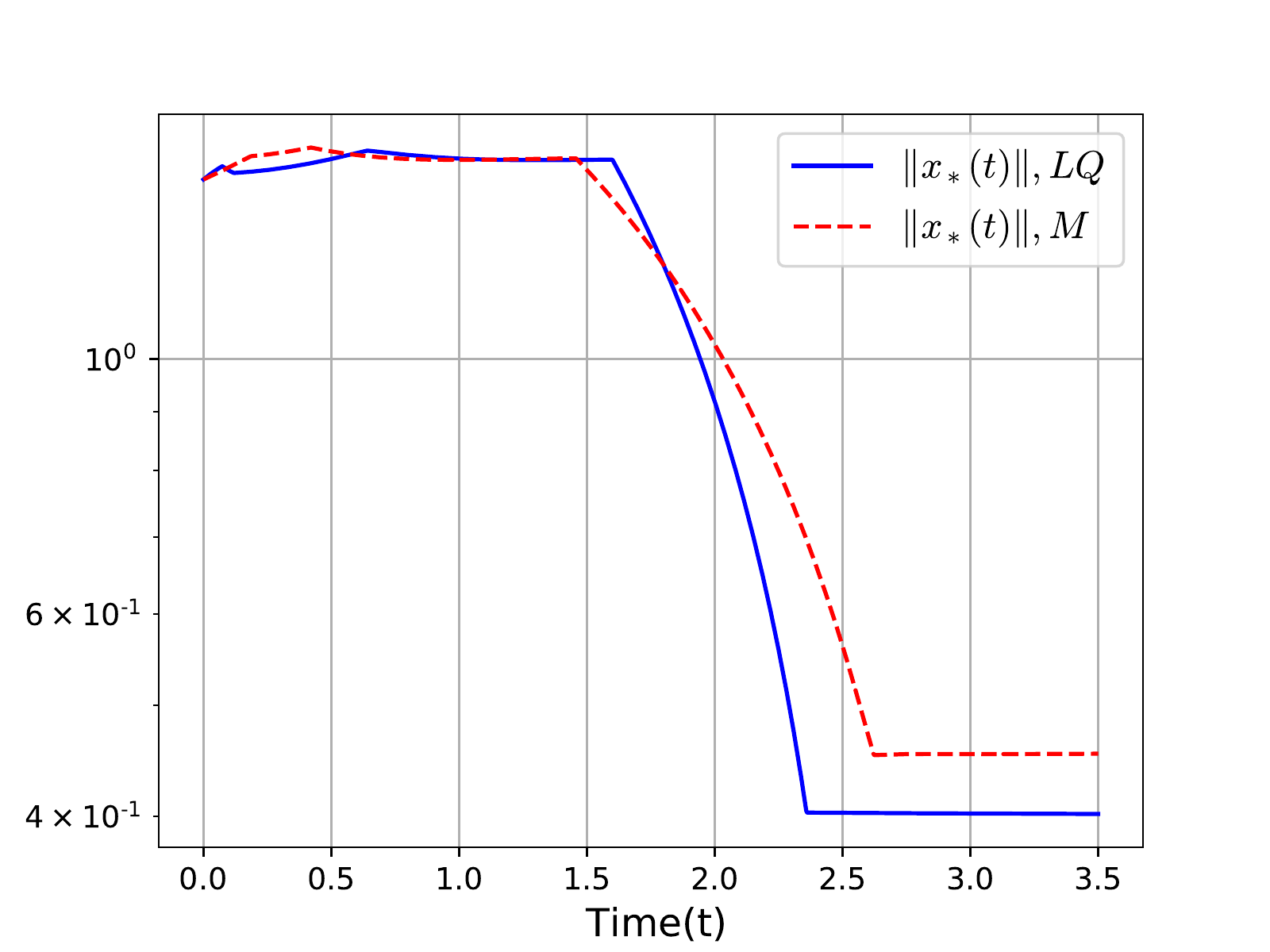}
			\caption{Norms of the states vs time for the LQ and the Mayer problems.}
			\label{fig:fig05}
		\end{figure}

%%%%%%%%%%%%%%%%%%%%%%%%%%%%%%%%%%%%%%%%%%%%%%%%%%%%%%%%%%%%%%%%%%%%%%%%%%%%%%%%
	\begin{example}
	\label{e:eg1}
		\textup{We first illustrate the optimal multiplexing for sparse LQ control problem described in \secref{s:LQ} (in \eqref{e:LQ problem}). Here we consider two linear systems, namely a harmonic oscillator (S1) and a linearised inverted pendulum on a moving cart (S2). The harmonic oscillator is described by the pair  $( A_1, B_1) = \left ( \pmat{0 & 1\\ -1 & 0 }, \pmat{ 0\\ 1}\right )$ and the inverted pendulum is described by the following fourth order model,}
	    \begin{equation}
			\label{eq6.1}
			(A_2, B_2) = \left( \pmat{0 & 1 & 0 & 0\\ 0 & 0 & \frac{-mg}{M} & 0\\ 0 & 0 & 0& 1\\ 0 & 0 &  \frac{(m+M)g}{LM} & 0 }, \pmat{ 0 \\ \frac{1}{M} \\ 0\\ \frac{-1}{LM}}\right).
        \end{equation}
		\textup{The parameter values are, \( m\)= 0.25 kg, \(M\)= 3 kg, \(g\)= 9.81 $m/s^2$, \(L\)= 2 m,  and weights, $ \lambda_1 =2, \lambda_2=1$, $ R_1 =2, R_2=2, Q_1 = 2I_2$, $ Q_2= \text{diag}((1,\, 5,\, 10,\, 10)\trnsp)$, $\hat{Q}_1 = 200I_2$, $ \hat{Q}_2= \text{diag}((10,\, 200,\, 200,\, 200)\trnsp)$, and the initial conditions, $x^1(0) =(1,0.5)\trnsp, x^2(0)=(0,\pi/10,0,0)\trnsp $ where superscripts indicate the corresponding system, $\hat{t}=3.5$ seconds. The convergence tolerance for all cases is kept at \(5\times 10^{-4}\).}

		\textup{The optimal control described by Theorem \ref{t:LQ} is applied in simulation to systems S1 and S2, and corresponding results shown in Figures \ref{fig:fig01}- \ref{fig:fig05}. These plots show comparison between results of examples \ref{e:eg1} and \ref{e:eg2}; curves corresponding to legend `LQ' are results of example (\ref{e:eg1}) considered here. Figure \ref{fig:fig01} shows multiplexed optimal controls \(t\mapsto u_*^1(t)\), \(t\mapsto u_*^2(t)\), where $u_*^1$  acts on the harmonic oscillator and $u_*^2$ acts on the inverted pendulum on a moving cart. It is evident from Figure \ref{fig:fig01} that at each time instant only one system is being controlled while the other evolves freely i.e., at time $t$ at most one of $\{u_*^1(t), u_*^2(t)\}$ is non-zero. Due to the additional sparsity requirement it can be seen that  $u_*^1$ and $u_*^2$ remain non zero only for a part of their activation times. For example, in the interval between \(0.15 s\) to \(0.65 s\) and from about \(2.4 s\) onwards, both the control commands are zero; this shows that the sparsity requirement is at work. A careful examination of both systems indicates sharp changes in solution trajectories whenever the control inputs $u_*^1$ and $u_*^2$ switch. Figure \ref{fig:fig02} shows the phase portrait of the harmonic oscillator (S1), wherein a higher terminal cost helps in bringing the states close to zero (the desired final position) with final values being $(x^1(\hat{t}),\dot{x}^1(\hat{t})) = (-0.0133, -0.0046)$.}

		\textup{Figure \ref{fig:fig03} plots the evolution of the inverted pendulum's angular position and angular velocity; both states remain close to zero for the entire duration of the simulation, which happens primarily because of the high running cost. $ \theta^2$ initially diverges from zero owing to the non-zero initial angular velocity but action of the controller ($u_*^2(t))$ aids recovery and both $ \theta^2(t)$ and $ \dot{\theta}^2(t)$ approach zero at final time. In this case as well, a higher terminal cost helps in bringing states close to zero, with final values being \((\theta^2(\hat{t}),\dot{\theta}^2(\hat{t})) = (1.60\times 10^{-4},  1.31\times10^{-3})\). In other words, the figures show approximate reachability of terminal condition (0,0). In Figure \ref{fig:fig04} we plot the position and velocity of the cart; while the velocity of the cart reaches zero around $1.6s$, the cart's position attains a non-zero value.}

		\textup{Figure \ref{fig:fig05} shows evolution of norm of all the six states (two for (S1) and four for (S2)). Initially the norm remains constant till about \(1.6 s\), and then decays rapidly to \(0.4\) units at \(2.3 s\) remaining constant beyond that. This coincides with the non-zero phase of $u_*^1$ acting on the harmonic oscillator. The norm attains a value of around \(0.4\) units, matching the cart's final position, while other states approach zero.}
	\end{example}

    \begin{example}
    	\label{e:eg2}
		\textup{In this example, we illustrate the optimal multiplexing for sparse Mayer problem described in \secref{s:Mayer} (described by \eqref{e:Mayer}). We consider the same system as in Example \ref{e:eg1} with identical parameters and initial conditions, and the terminal values of all states are kept at \(0\). The parameters $\lambda_1,\, \lambda_2,\, \hat{Q}_1 ,\, \hat{Q}_2$ and final time $\hat{t}$ are also the same as for Example \ref{e:eg1} with $\subadmact{1} = [-1,1]$ and $\subadmact{2} = [-1,1]$. The convergence tolerance in all the experiments is \(5\times 10^{-4}\).}

		\textup{The optimal control described by Theorem \ref{t:Mayer} is applied for simulation of the aforementioned system and corresponding results are shown in Figures \ref{fig:fig01} - \ref{fig:fig05}, with the curves corresponding to legend ``M'' depicting the results of the current example. From these plots, it is evident that the evolution of all states follows a pattern similar to that of the LQ case. The significant differences are due to the absence of a running quadratic cost (\(Q=0\) and \(R=0\)) and the presence of bounded admissible action sets.}

		\textup{We observe that due to bounds on the the control actions, the optimal controls $u_*^1$ and $u_*^2$ switchd between \(\{-1\}\), \(\{0\}\) and \(\{1\}\), i.e., both of them have a \textit{bang-off-bang} profile as expected. Also, in this case, the controllers are active over a longer time span as compared to the LQ problem, and the target achievement is slightly inferior than that of the LQ optimal controller. However, increasing the weight $R$ on the control in the LQ cost leads to similar controller performance in both cases. As observed in Example \ref{e:eg1}, sharp changes in trajectories of both the systems are clearly visible whenever $u_*^1$ or $u_*^2$ switch.}

		\textup{In this case as well, the norm shows a continuous decay starting before $t=1.5 s$ till $t=2.6s$ and reaches around \(0.45\) units which coincides with the activation of $u_*^1$ acting on the harmonic oscillator. The norm remains constant beyond  $t=2.6s$. From Figure \ref{fig:fig01} it is clear that at each instant only one system is being controlled. A comparison of Figure \ref{fig:fig03} and Figure \ref{fig:fig04} indicates that the LQ optimal controller performs better in the context of reaching the desired final condition. However, neither of these examples solves the exact reachability problem~\eqref{e:reach problem}.}
	\end{example}                               

\end{multicols}

\begin{multicols}{2}

%%%%%%%%%%%%%%%%%%%%%%%%%%%%%%%%%%%%%%%%%%%%%%%%%%%%%%%%%%%%%%%%%%%%%%%%%%%%%%%%
	\bibliographystyle{siam}
	\bibliography{refr}
%%%%%%%%%%%%%%%%%%%%%%%%%%%%%%%%%%%%%%%%%%%%%%%%%%%%%%%%%%%%%%%%%%%%%%%%%%%%%%%%

\appendix

%%%%%%%%%%%%%%%%%%%%%%%%%%%%%%%%%%%%%%%%%%%%%%%%%%%%%%%%%%%%%%%%%%%%%%%%%%%%%%%%
	\section{A nonsmooth Pontryagin maximum principle}
	\label{s:Clarke extended PMP}
%%%%%%%%%%%%%%%%%%%%%%%%%%%%%%%%%%%%%%%%%%%%%%%%%%%%%%%%%%%%%%%%%%%%%%%%%%%%%%%%
		We need the following adaptation of \cite[Theorem 22.26]{ref:Cla-13}:\footnote{The assumptions of \cite[Theorem 22.26]{ref:Cla-13} are considerably weaker than those of Theorem \ref{t:Clarke extended PMP}, but Theorem \ref{t:Clarke extended PMP} is sufficient for our purposes here.}
		\begin{theorem}
		\label{t:Clarke extended PMP}
			Let \(\tfin\in\:]0, +\infty[\) and let \(\admact\subset\R^m\) be a non-empty Borel measurable set. Let a lower semicontinuous instantaneous cost function \(\R^d\times\admact\ni(\xi, \mu)\mapsto\Lambda(\xi, \mu)\in\R\), with \(\Lambda\) continuously differentiable in \(\xi\) for every fixed \(\mu\),\footnote{Recall that a map \(\varphi:X\lra\R\) from a topological space \(X\) into the real numbers is said to be lower semicontinuous if for every \(c\in\R\) the set \(\{z\in\R^d\mid \varphi(z) \le c\}\) is closed, and a map \(\psi:X\lra\R\) is said to be upper semicontinuous if \(-\psi\) is lower semicontinuous.} and a continuously differentiable terminal cost function \(\ell:\R^d\times\R^d\lra\R\) be given. Consider the optimal control problem
			\begin{equation}
			\label{e:Clarke extended opt problem}
			\begin{aligned}
				& \minimize_{u}	&& \quad  \ell\bigl(x(\tinit), x(\tfin)\bigr) + \int_{\tinit}^{\tfin} \Lambda\bigl(x(t), u(t)\bigr) \, \dd t \\
				& \sbjto		&& \quad  \begin{cases}
						\dot x(t) = f\bigl( x(t), u(t)\bigr)\quad \text{for a.e.\ }t\in[\tinit, \tfin],\\
						u(t) \in\admact\quad\text{for a.e.\ }t \in [\tinit, \tfin],\\
						u\text{ Lebesgue measurable},\\
						\bigl(x(\tinit), x(\tfin)\bigr) \in E \subset\R^d\times\R^d,
					\end{cases}
			\end{aligned}
			\end{equation}
			where \(f:\R^d\times\R^m\lra\R^d\) is continuously differentiable, and \(E\) is a closed set. For a real number \(\eta\), we define the \emph{Hamiltonian} $H^\eta$ by 
			\[
				H^\eta(p, x, u) \Let \inprod{p}{f(x, u)} - \eta \Lambda(x, u),\quad (p, x, u)\in\R^d\times\R^d\times\admact.
			\]
			If \([\tinit, \tfin]\ni t\mapsto \bigl( x\opt(t), u\opt(t) \bigr)\in\R^d\times\admact\) is a local minimizer of \eqref{e:Clarke extended opt problem}, then there exist an absolutely continuous map \( p: [\tinit, \tfin] \lra \R^d \) and a scalar \( \eta \) equal to \(0\) or \(1\), satisfying the \emph{nontriviality condition}
			\begin{equation}
			\label{e:Clarke:nontriviality}
				\bigl(\eta, p(t)\bigr) \neq 0 \quad \text{for all } t \in [\tinit, \tfin],
			\end{equation}
			the \emph{transversality condition}%\todo[inline, color=green!40]{check out the precise meanings}
			\begin{equation}
			\label{e:Clarke:transversality}
				\bigl(p(\tinit), - p(\tfin)\bigr) \in \eta \nabla\ell\bigl( x\opt(\tinit), x\opt(\tfin)\bigr) + N_E^L \bigl(x\opt(\tinit), x\opt(\tfin)\bigr),
			\end{equation}
			where \(\nabla\ell\) is the gradient of \(\ell\) and \(N_E^L\bigl(x\opt(\tinit), x\opt(\tfin)\bigr)\) is the limiting normal cone to \(E\) at the point \(\bigl(x\opt(\tinit), x\opt(\tfin)\bigr)\),\footnote{The \emph{limiting normal cone} to a closed subset \(S\) of \(\R^\nu\) is defined by means of a topological closure operation applied to the \emph{proximal normal cone} to the set \(S\); see, e.g., \cite[p.\ 240]{ref:Cla-13} for the definition of the proximal normal cone, and \cite[p.\ 244]{ref:Cla-13} for the definition of the limiting normal cone.} the \emph{adjoint equation}
			\begin{equation}
			\label{e:Clarke:adjoint}
				- \dot p(t) = \partial_x H^\eta \bigl( p(t), \boldsymbol\cdot, u\opt(t) \bigr) (x\opt(t))\quad\text{for a.e.\ }t\in[\tinit, \tfin],
			\end{equation}
			the \emph{Hamiltonian maximum condition}
			\begin{equation}
			\label{e:Clarke:Hamiltonian max}
				H^\eta\bigl( p(t), x\opt(t), u\opt(t) \bigr) = \sup_{v\in\admact} H^\eta\bigl( p(t), x\opt(t), v \bigr)\quad\text{for a.e.\ }t\in[\tinit, \tfin],
			\end{equation} 
			as well as  the \emph{constancy of the Hamiltonian}
			\begin{equation}
			\label{e:Clarke:Hamiltonian constancy}
				H^\eta\bigl(p(t), x\opt(t), u\opt(t)\bigr) = h\quad\text{for a.e.\ }t\in[\tinit, \tfin] \text{ and some }h\in\R.
			\end{equation}
		\end{theorem}

		The quadruple \(\bigl(\eta, p, x\opt, u\opt\bigr)\) is known as the \emph{extremal lift} of the optimal state-action trajectory \([\tinit, \tfin]\ni t\mapsto \bigl(x\opt(t), u\opt(t)\bigr)\). The number \(\eta\) is called the \emph{abnormal multiplier}. The \emph{abnormal} case --- when \(\eta = 0\) --- may arise, e.g., when the constraints of the optimal control problem are so tight that the cost function plays no r\^ole in determining the solution. For instance, we have an abnormal case when the optimal solution \(t\mapsto \bigl(x\opt(t), u\opt(t)\bigr)\) is ``isolated'' in the sense that there is no other solution satisfying the end-point constraints in the vicinity --- as measured by the supremum norm --- of the optimal solution.

		Theorem \ref{t:Clarke extended PMP} will yield characterizations of the controls that solve the optimal control problems \eqref{e:reach problem} and \eqref{e:LQ problem}; it constitutes the backbones of the proofs below.

%%%%%%%%%%%%%%%%%%%%%%%%%%%%%%%%%%%%%%%%%%%%%%%%%%%%%%%%%%%%%%%%%%%%%%%%%%%%%%%%
	\section{Multiplexed sparsest reachability: proofs}
	\label{s:proof:reach}
%%%%%%%%%%%%%%%%%%%%%%%%%%%%%%%%%%%%%%%%%%%%%%%%%%%%%%%%%%%%%%%%%%%%%%%%%%%%%%%%
		Recall that in \secref{s:reach} we defined \(\lambda_k' = \sum_{\ell\in\{1, \ldots, N\}\setmin\{k\}} \lambda_\ell\) for each \(k = 1, \ldots, N\), and \(\wt\lambda = \sum_{k=1}^N \lambda_k\).

		\textbf{Proof of Theorem \ref{t:reach}}
			Notice first that Theorem \ref{t:Clarke extended PMP} applies directly to the optimal control problem \eqref{e:reach problem}. Indeed, 
			\begin{itemize}[label=\(\circ\), leftmargin=*]
				\item \(\admactreach\) is a finite union of compact sets, and is therefore Borel measurable;
				\item the dynamics \(f\) is given by the linear control system \eqref{e:joint system} and is therefore smooth;
				\item the instantaneous cost function \(\Lambda(\xi, \mu) = - \sum_{k=1}^N \lambda_k \indic{\{0\}}(\mu^k)\) is independent of the space variable \(\xi\) and is lower semicontinuous in \(\mu\);
				\item the terminal cost \(\ell\) is identically \(0\); and
				\item the boundary constraint set \(E = \{(\stinit, \stfin)\}\) is a singleton, and is therefore closed.
			\end{itemize}

			Let the state-action trajectory \([\tinit, \tfin]\ni t\mapsto \bigl(\st\opt(t), \con\opt(t)\bigr)\) be a local minimizer of \eqref{e:reach problem}. For a real number \(\eta\) we define the Hamiltonian
			\begin{equation}
			\label{e:reach:Hamiltonian}
			\begin{aligned}
				& (\R^{d_1}\times\cdots\times\R^{d_N})\times(\R^{d_1}\times\cdots\times\R^{d_N})\times\admactreach\\
				& \quad \ni \left(\pmat{\subadj{1}\\ \smash[t]{\vdots} \\ \subadj{N}}, \pmat{\xi^1\\\smash[t]{\vdots}\\\xi^N}, \pmat{\mu^1\\\smash[t]{\vdots}\\\mu^N}\right) \teL (\adj, \xi, \mu) \mapsto\\
				& \quad \Ham^\eta(\adj, \xi, \mu) \Let \inprod{\adj}{A\xi + B\mu} + \eta \sum_{k=1}^N \lambda_k \indic{\{0\}}(\mu^k)\\
				& = \sum_{k=1}^N \left( \inprod{\subadj{k}}{A_k \xi^k + B_k \mu^k} + \eta \lambda_k \indic{\{0\}}(\mu^k) \right) \in\R,
			\end{aligned}
			\end{equation}
			where we have employed the block-diagonal structure of \(A\) and \(B\) to arrive at the last equality.

			By Theorem \ref{t:Clarke extended PMP}, there exists an absolutely continuous map \([\tinit, \tfin]\ni t\mapsto \adj(t)\in\R^d\) (called an \emph{adjoint trajectory},) that, in view of the adjoint equation \eqref{e:Clarke:adjoint}, solves 
			\[
				-\dot \adj(t) = \partial_\xi \Ham^\eta\bigl(\adj(t), \cdot, \con\opt(t)\bigr)\bigl(\st\opt(t)\bigr) = A\trnsp \adj(t)\quad \text{for a.e.\ }t\in[\tinit, \tfin],
			\]
			or, in terms of the individual \(N\) maps \((\subadj{k})_{k=1}^N\) obtained by projecting \(\adj\) at each time to appropriate factors in an obvious way,
			\[
				-\tdsubadj{k}(t) = A_k\trnsp\subadj{k}(t)\quad \text{for a.e.\ }t\in[\tinit, \tfin],\quad k = 1, \ldots, N.
			\]
			Note that linearity of the right-hand sides ensure that there is a unique adjoint trajectory. The boundary conditions for the adjoint \(\adj\) are obtained from the transversality conditions \eqref{e:Clarke:transversality}, and in our problem they turn out to be
			\[
				\bigl(\adj(\tinit), -\adj(\tfin)\bigr) \in \R^d\times\R^d,
			\]
			in other words, the boundary conditions of \(\adj\) are free.

			The optimal control actions as functions of time are such that they satisfy the Hamiltonian maximum condition \eqref{e:Clarke:Hamiltonian max}. For our problem this condition is: for a.e.\ \(t\in[\tinit, \tfin]\),
			\begin{align*}
				& \con\opt(t) \in \argmax_{\mu\in\admactreach}\Ham^\eta\bigl(\adj(t), \st\opt(t), \mu\bigr)\\
				& \quad = \argmax_{\mu\in\admactreach} \Biggl\{ \sum_{k=1}^N \left( \inprod{\subadj{k}(t)}{A_k \subst{k}\opt(t) + B_k \mu^k} + \eta \lambda_k \indic{\{0\}}(\mu^k) \right) \Biggr\}\\
				& \quad = \argmax_{\mu\in\admactreach} \Biggl\{ \sum_{k=1}^N \left( \inprod{B_k\trnsp\subadj{k}(t)}{\mu^k} + \eta \lambda_k \indic{\{0\}}(\mu^k) \right) \Biggr\}.
			\end{align*}
			Denoting by \(S\subset[\tinit, \tfin]\) the full-measure set on which the preceding membership of \(\con\opt\) holds, we fix \(t\in S\). For this \(t\), any arguemnt \(\mu \Let (\mu^k)_{k=1}^N\in\admactreach\) of the map
			\[
				\admactreach\ni\pmat{\mu^1\\\smash[t]{\vdots}\\\mu^N}\mapsto \sum_{k=1}^N \biggl( \inprod{B_k\trnsp \adj{k}(t)}{\mu^k} + \eta\lambda_k \indic{\{0\}}(\mu^k) \biggr) \in \R
			\]
			has at most one non-zero \(\mu^k\) due to the ``star-shaped'' structure of the set \(\admactreach\) defined in \eqref{e:admact reach}. For each \(k = 1, \ldots, N\), we let
			\begin{multline*}
				\subadmact{k}\ni v\mapsto \phi^\eta_k(t, v)	\Let \inprod{B_k\trnsp \subadj{k}(t)}{v}\\
					+ \eta\bigl( \lambda_k \indic{\{0\}}(v) + \lambda_k' \bigr) \in \R,
			\end{multline*}
			and note that \(\phi^\eta_k(t, \cdot)\) is upper semicontinuous; since the set \(\subadmact{k}\) is compact by assumption,
			\[
				\ol\phi^\eta_k(t) \Let \sup_{v\in\subadmact{k}} \phi^\eta_k(t, v)
			\]
			is attained on \(\subadmact{k}\) by Weierstrass's theorem \cite[Exercise 2.14]{ref:Cla-13}. We let \(\Phi^\eta(t, k)\) denote the non-empty set of maximizers of \(\phi^\eta_k(t, \cdot)\), \(k = 1, \ldots, N\); i.e.,
			\begin{equation}
			\label{e:Phieta def}
				\Phi^\eta(t, k) \Let \argmax_{v\in\subadmact{k}} \phi^\eta_k(t, v),\quad k = 1, \ldots, N.
			\end{equation}
			Informally, at the time \(t\) fixed above, we get the finite sequence \(\bigl(\ol\phi^\eta_k(t)\bigr)_{k=1}^N\) of real numbers, and this finite sequence has a maximum element, say \(\ol\phi^\eta_{k\opt}(t)\); the optimal control action \(\con\opt(t) = \bigl(\subcon{k}\opt(t)\bigr)_{k=1}^N\) must be such that \(\subcon{k\opt}\opt(t) \in \Phi^\eta(t, k\opt)\) and \(\subcon{\ell}\opt(t) = 0\) for all \(\ell\in\{1, \ldots, N\}\setmin\{k\opt\}\). By letting \(t\) range over \(S\), we have a characterization of \(u\opt\) on \(S\). The behaviour of \(u\opt\) on \([\tinit, \tfin]\setmin S\) can, of course, be arbitrary. Formally, defining
			\begin{equation}
			\label{e:reachmuxset def}
				\reachmuxset_\eta(t) \Let \begin{dcases*}
					\argmax_{k\in\{1, \ldots, N\}} \ol\phi^\eta_k(t)	& for \(t\in S\),\\
					\{1, \ldots, N\}	& for \(t\in[\tinit, \tfin]\setmin S\),
				\end{dcases*}
			\end{equation}
			we arrive at a family of non-empty subsets of \(\{1, \ldots, N\}\) parametrized by \(t\in[\tinit, \tfin]\); in other words, \(\reachmuxset_\eta\) is a set-valued map from \([\tinit, \tfin]\) into the power set of \(\{1, \ldots, N\}\). Given \(\reachmuxset_\eta\), any map (commonly known as a \emph{selector} of the set-valued map \(\reachmuxset_\eta\),)
			\[
				[\tinit, \tfin]\ni t\mapsto \reachmux_\eta(t)\in\reachmuxset_\eta(t)
			\]
			gives us an admissible multiplexer. (Note that the Axiom of Choice \cite[p.\ 8]{ref:DiB-02}, guarantees that there always exists such a selector, and therefore, a multiplexer.) It follows that the set \((\subcon{k}\opt)_{k=1}^N\) of optimal controls constituting \(\con\opt\) satisfies
			\begin{equation}
			\label{e:reach muxcon}
			\begin{aligned}
				& \subcon{k}\opt(t)\in \begin{cases}
						\Phi^\eta\bigl(t, \reachmux_\eta(t)\bigr)	& \text{if }k = \reachmux_\eta(t),\\
						\{0\}										& \text{otherwise},
					\end{cases}
					\quad k = 1, \ldots, N,\\
				& \text{for all }t\in[\tinit, \tfin].
			\end{aligned}
			\end{equation}

			It is time to provide more precise descriptions of the sets \(\Phi^\eta(t, k)\) defined in \eqref{e:Phieta def}. As asserted by Theorem \ref{t:Clarke extended PMP}, only the two cases of \(\eta = 0\) or \(\eta = 1\) arise. If \(\eta = 0\), then for each \(t\in S\)
			\[
				\Phi^0(t, k) = \argmax_{v\in\subadmact{k}} \inprod{B_k\trnsp \subadj{k}(t)}{v}\quad \text{for }k = 1, \ldots, N,
			\]
			and
			\[
				\reachmux(t) \in \argmax_{k\in\{1, \ldots, N\}} \max_{v\in\subadmact{k}} \inprod{B_k\trnsp \subadj{k}(t)}{v}.
			\]
			If \(\eta = 1\), then for each \(t\in S\)
			\begin{align*}
				& \Phi^1(t, k) = \argmax_{v\in\subadmact{k}}\left\{\inprod{B_k\trnsp \subadj{k}(t)}{v} + \lambda_k\indic{\{0\}}(v)\right\}\\
				& = \begin{dcases}
					\argmax_{v\in\subadmact{k}}\inprod{B_k\trnsp\subadj{k}(t)}{v}	& \text{if }\max_{v\in\subadmact{k}} \inprod{B_k\trnsp\subadj{k}(t)}{v} > \lambda_k,\\
					\{0\}\cup\argmax_{v\in\subadmact{k}}\inprod{B_k\trnsp\subadj{k}(t)}{v}	& \text{if }\max_{v\in\subadmact{k}} \inprod{B_k\trnsp\subadj{k}(t)}{v} = \lambda_k,\\
					\{0\}	& \text{otherwise}.
				\end{dcases}
				\\
				& \quad \text{for }k = 1, \ldots, N,
			\end{align*}
			and
			\[
				\reachmux(t)\in \argmax_{k \in\{1, \ldots, N\}} \max_{v\in\subadmact{k}} \Bigl\{ \inprod{B_k\trnsp\subadj{k}(t)}{v} + \lambda_k \indic{\{0\}}(v) + \lambda_k' \Bigr\}.
			\]
			Note that \(\lambda_k'\) is a constant in the definition of \(\phi^1_k\), and plays no r\^ole in the determination of the set \(\Phi^1(t, k)\). The value of \(\ol\phi^1_k(t)\) at each \(t\in S\), however, depends on this constant, and therefore, so does the set-valued map \(t\mapsto \reachmuxset_1(t)\) defined in \eqref{e:reachmuxset def}, and therefore, also the map \(t\mapsto\reachmux_1(t)\). Moreover, since \(\con\opt\) is measurable, so is \(\reachmux\).

			The constancy of the Hamiltonian \eqref{e:Clarke:Hamiltonian constancy} gives the final assertion of the theorem, thereby completing the proof.\hfill{}\(\square\)

		\textbf{Proof of Corollary \ref{c:reach scalar}}
			We retain the notations introduced in the proof of Theorem \ref{t:reach} above, and note that the only details that need to be supplied here are the sets \(\Phi^0(t, k)\) and \(\Phi^1(t, k)\) for each \(t\in S\) and \(k = 1, \ldots, N\), corresponding to \(\eta = 0\) and \(\eta = 1\), respectively. To that end, note that if \(\eta = 0\), then
			\begin{align*}
				\Phi^0(t, k)	& = \argmax_{v\in[-1, 1]} \phi^0_k(t, v) = \argmax_{v\in[-1, 1]} B_k\trnsp \subadj{k}(t) v\\
					& = \sgn\big(B_k\trnsp \subadj{k}(t)\bigr),
			\end{align*}
			where the set-valued map \(\sgn\) defined in \eqref{e:sgn def}. If \(\eta = 1\), then
			\begin{align*}
				\Phi^1(t, k)	& = \argmax_{v\in[-1, 1]} \phi^1_k(t, v)\\
					& = \argmax_{v\in[-1, 1]}\Bigl\{ B_k\trnsp\subadj{k}(t) v + \lambda_k \indic{\{0\}}(v) + \lambda_k' \Bigr\}\\
					& = \begin{dcases}
						\sgn\bigl(B_k\trnsp \subadj{k}(t)\bigr)	& \text{if }\abs{B_k\trnsp\subadj{k}(t)} > \lambda_k,\\
						\{0\}\cup\sgn\bigl(B_k\trnsp \subadj{k}(t)\bigr)	& \text{if }\abs{B_k\trnsp\subadj{k}(t)} = \lambda_k,\\
						\{0\}	& \text{otherwise}.
					\end{dcases}
			\end{align*}
			The assertion follows at once.\hfill{}\(\square\)

%%%%%%%%%%%%%%%%%%%%%%%%%%%%%%%%%%%%%%%%%%%%%%%%%%%%%%%%%%%%%%%%%%%%%%%%%%%%%%%%
	\section{Multiplexed sparse LQ control: proof}
	\label{s:proof:LQ}
%%%%%%%%%%%%%%%%%%%%%%%%%%%%%%%%%%%%%%%%%%%%%%%%%%%%%%%%%%%%%%%%%%%%%%%%%%%%%%%%
		Recall that in \secref{s:LQ} we defined \(\lambda_k' = \sum_{\ell\in\{1, \ldots, N\}\setmin\{k\}} \lambda_\ell\) for each \(k = 1, \ldots, N\), and \(\wt\lambda = \sum_{k=1}^N \lambda_k\).

		\textbf{Proof of Theorem \ref{t:LQ}}
			Notice first that Theorem \ref{t:Clarke extended PMP} applies directly to the optimal control problem \eqref{e:LQ problem}; the details being similar those elaborated in the Proof of Theorem \ref{t:reach} in \secref{s:proof:reach}, we omit them in the interest of brevity. Let the state-action trajectory \([\tinit, \tfin]\ni t\mapsto \bigl(\st\opt(t), \con\opt(t)\bigr)\) be a local minimizer of \eqref{e:reach problem}. For a real number \(\eta\) we define the Hamiltonian
			\begin{equation}
			\label{e:LQ:Hamiltonian}
			\begin{aligned}
				& (\R^{d_1}\times\cdots\times\R^{d_N})\times(\R^{d_1}\times\cdots\times\R^{d_N})\times\admactLQ\\
				& \quad \ni \left(\pmat{\subadj{1}\\ \smash[t]{\vdots} \\ \subadj{N}}, \pmat{\xi^1\\\smash[t]{\vdots}\\\xi^N}, \pmat{\mu^1\\\smash[t]{\vdots}\\\mu^N}\right) \teL (\adj, \xi, \mu) \mapsto\\
				& \quad \Ham^\eta(\adj, \xi, \mu) \Let \inprod{\adj}{A\xi + B\mu}\\
				& \quad - \eta\biggl( \tfrac{1}{2}\Bigl(\inprod{\xi}{Q\xi} + \inprod{\mu}{R\mu}\Bigr) - \sum_{k=1}^N \lambda_k \indic{\{0\}}(\mu^k)\biggr)\\
				& = \sum_{k=1}^N \biggl( \inprod{\subadj{k}}{A_k \xi^k + B_k \mu^k}\\
				& \qquad - \eta\Bigl( \tfrac{1}{2}\Bigl(\inprod{\xi^k}{Q_k \xi^k} + \inprod{\mu^k}{R_k \mu^k}\Bigr) - \lambda_k \indic{\{0\}}(\mu^k) \Bigr) \biggr) \in\R,
			\end{aligned}
			\end{equation}
			where we have employed the block-diagonal structure of \(A\) and \(B\) to arrive at the last equality.

			By Theorem \ref{t:Clarke extended PMP}, there exists an absolutely continuous map \([\tinit, \tfin]\ni t\mapsto \adj(t)\in\R^d\) (called the \emph{adjoint trajectory},) that, in view of the adjoint equation \eqref{e:Clarke:adjoint}, solves 
			\begin{equation}
			\label{e:LQ adj}
			\begin{aligned}
				-\dot \adj(t)	& = \partial_\xi \Ham^\eta\bigl(\adj(t), \cdot, \con\opt(t)\bigr)\bigl(\st\opt(t)\bigr)\\
					& = A\trnsp \adj(t) + \eta Q \st\opt(t)
			\end{aligned}
			\quad \text{for a.e.\ }t\in[\tinit, \tfin],
			\end{equation}
			or, in terms of the individual \(N\) maps \((\subadj{k})_{k=1}^N\) obtained by projecting \(\adj\) at each time to appropriate factors in an obvious way,
			\[
				-\tdsubadj{k}(t) = A_k\trnsp\subadj{k}(t) + \eta Q_k \subst{k}\opt(t)\quad \text{for a.e.\ }t\in[\tinit, \tfin],\quad k = 1, \ldots, N.
			\]
			Linearity of the right-hand sides ensure that there is a unique adjoint trajectory. The boundary conditions for the adjoint \(\adj\) are obtained from the transversality conditions \eqref{e:Clarke:transversality}, and in our problem they turn out to be
			\[
				\bigl(\adj(\tinit), -\adj(\tfin)\bigr) \in \R^d\times\{\eta \hat Q \st\opt(\tfin)\};
			\]
			in other words, the initial condition of \(\adj\) is free, and the final condition is \(\adj(\tfin) = -\eta \hat Q\st\opt(\tfin)\).

			Theorem \ref{t:Clarke extended PMP} admits only two cases of \(\eta\) --- \(0\) or \(1\). We claim that the case of \(\eta = 0\) does not arise in \eqref{e:LQ problem}. Indeed, if \(\eta = 0\), then the final boundary condition of each \(\subadj{k}\) is \(0\), leading to \(\adj(\tfin) = 0\), and the forcing term on the right-hand side of \eqref{e:LQ adj} also vanishes. In view of the resulting linearity of \eqref{e:LQ adj} with final condition equal to \(0\), the entire trajectory \(\adj\) vanishes. This contradicts the non-triviality condition \eqref{e:Clarke:nontriviality} of Theorem \ref{t:Clarke extended PMP}. Therefore, \(\eta\) must be equal to \(1\), to which we commit and henceforth write \(\Ham\) instead of \(\Ham^1\).

			The optimal control actions as functions of time are such that they satisfy the Hamiltonian maximum condition \eqref{e:Clarke:Hamiltonian max}. For our problem this condition is: for a.e.\ \(t\in[\tinit, \tfin]\),
			\begin{align*}
				& \con\opt(t) \in \argmax_{\mu\in\admactLQ}\Ham\bigl(\adj(t), \st\opt(t), \mu\bigr)\\
				& \quad = \argmax_{\mu\in\admactLQ} \Biggl\{\sum_{k=1}^N \biggl( \inprod{\subadj{k}(t)}{A_k \subst{k}\opt(t) + B_k \mu^k}\\
				& \qquad\qquad\qquad - \Bigl( \tfrac{1}{2}\inprod{\xi^k}{Q_k \xi^k} + \tfrac{1}{2}\inprod{\mu^k}{R_k \mu^k} - \lambda_k \indic{\{0\}}(\mu^k) \Bigr) \biggr) \Biggl\}\\
				& \quad = \argmax_{\mu\in\admactLQ} \Biggl\{ \sum_{k=1}^N \biggl( \inprod{B_k\trnsp\subadj{k}(t)}{\mu^k}- \tfrac{1}{2}\inprod{\mu^k}{R_k \mu^k} + \lambda_k \indic{\{0\}}(\mu^k) \biggr) \Biggr\}.
			\end{align*}
			Denoting by \(S\subset[\tinit, \tfin]\) the full-measure set on which the preceding membership of \(\con\opt\) holds, we fix \(t\in S\). For this \(t\), any argument \(\mu \Let (\mu^k)_{k=1}^N\in\admactLQ\) of the map
			\begin{align*}
				\admactLQ\ni\pmat{\mu^1\\\smash[t]{\vdots}\\\mu^N}\mapsto \sum_{k=1}^N & \biggl( \inprod{B_k\trnsp\subadj{k}(t)}{\mu^k}\\
					& - \tfrac{1}{2} \inprod{\mu^k}{R_k \mu^k} + \lambda_k \indic{\{0\}}(\mu^k) \biggr) \in\R
			\end{align*}
			has at most one non-zero \(\mu^k\) due to the star-shaped conical structure of the set \(\admactLQ\) defined in \eqref{e:admact LQ}. For each \(k = 1, \ldots, N\), we let
			\begin{multline*}
				\R^{m_k}\ni v\mapsto \psi_k(t, v) \Let \inprod{B_k\trnsp \subadj{k}(t)}{v} - \tfrac{1}{2}\inprod{v}{R_k v}\\
					+ \lambda_k \indic{\{0\}}(v) + \lambda_k' \in \R.
			\end{multline*}
			Note that the first two terms comprising \(\psi_k(t, \cdot)\) define a strictly concave function due to positive definiteness of the matrix \(R_k\), and the second two terms are bounded above. Moreover, \(\psi_k(t, \cdot)\) is upper semicontinuous, and satisfies \(\psi_k(t, v)\xrightarrow[\norm{v}\to+\infty]{} -\infty\); therefore,
			\begin{equation}
			\label{e:olpsik def}
				\ol\psi_k(t) \Let \sup_{v\in\R^{m_k}} \psi_k(t, v)
			\end{equation}
			is attained on \(\R^{m_k}\) by standard arguments that ultimately rely on Weierstrass's theorem \cite[Exercise 2.14]{ref:Cla-13}. We let \(\Psi(t, k)\) denote the non-empty set of maximizers of \(\psi_k(t, \cdot)\), \(k = 1, \ldots, N\); i.e.,
			\begin{equation}
			\label{e:Psi def}
				\Psi(t, k) \Let \argmax_{v\in\R^{m_k}} \psi_k(t, v),\quad k = 1, \ldots, N.
			\end{equation}
			As we shall see momentarily, more than one maximizers of each \(\psi_k(t, \cdot)\) may exist. Informally, at the time \(t\) fixed above, we get the finite sequence \(\bigl(\ol\psi_k(t)\bigr)_{k=1}^N\) of real numbers, and this finite sequence has a maximum element, say \(\ol\psi_{k\opt}(t)\); the optimal control action \(\con\opt(t) = \bigl(\subcon{k}\opt(t)\bigr)_{k=1}^N\) must be such that \(\subcon{k\opt}\opt(t) \in \Psi(t, k\opt)\) and \(\subcon{\ell}\opt(t) = 0\) for all \(\ell\in\{1, \ldots, N\}\setmin\{k\opt\}\). By letting \(t\) range over \(S\), we have a characterization of \(u\opt\) on \(S\). The behaviour of \(u\opt\) on \([\tinit, \tfin]\setmin S\) can, of course, be arbitrary. Formally, defining
			\begin{equation}
			\label{e:LQmuxset def}
				\LQmuxset(t) \Let \begin{dcases}
					\argmax_{k\in\{1, \ldots, N\}} \ol\psi_k(t)	& \text{for }t\in S,\\
					\{1, \ldots, N\}	& \text{for }t\in[\tinit, \tfin]\setmin S,
				\end{dcases}
			\end{equation}
			we arrive at a family of non-empty subsets of \(\{1, \ldots, N\}\) parametrized by \(t\in[\tinit, \tfin]\); in other words, \(\LQmuxset\) is a set-valued map from \([\tinit, \tfin]\) into the power set of \(\{1, \ldots, N\}\). Given \(\LQmuxset\), any map (i.e., any \emph{selector} of the set-valued map \(\LQmuxset\),)
			\[
				[\tinit, \tfin]\ni t\mapsto \LQmux(t)\in\LQmuxset(t)
			\]
			gives us an admissible multiplexer. (The Axiom of Choice \cite[p.\ 8]{ref:DiB-02} guarantees the existence of such a selector, and therefore, a multiplexer.) It follows that the set \((\subcon{k}\opt)_{k=1}^N\) of optimal controls constituting \(\con\opt\) satisfies
			\begin{equation}
			\label{e:LQ muxcon}
			\begin{aligned}
				& \subcon{k}\opt(t)\in \begin{cases}
						\Psi\bigl(t, \LQmux(t)\bigr)	& \text{if }k = \LQmux(t),\\
						\{0\}									& \text{otherwise},
					\end{cases}
					\quad k = 1, \ldots, N,\\
				& \text{for all }t\in[\tinit, \tfin].
			\end{aligned}
			\end{equation}

			We provide more precise descriptions of the sets \(\Psi(t, k)\) defined in \eqref{e:Psi def}: For each \(t\in S\) and \(k = 1, \ldots, N\),
			\begin{align*}
				\Psi(t, k)	& = \argmax_{v\in\R^{m_k}}\left\{\inprod{B_k\trnsp \subadj{k}(t)}{v} - \tfrac{1}{2}\inprod{v}{R_k v} + \lambda_k\indic{\{0\}}(v)\right\}\\
					& = \begin{cases}
						\bigl\{R_k\inverse B_k\trnsp \subadj{k}(t)\bigr\}			& \text{if } \norm{\subadj{k}(t)}_{B_k R_k\inverse B_k\trnsp}^2 > 2\lambda_k,\\
						\{0\}\cup\bigl\{R_k\inverse B_k\trnsp \subadj{k}(t)\bigr\}	& \text{if } \norm{\subadj{k}(t)}_{B_k R_k\inverse B_k\trnsp}^2 = 2\lambda_k,\\
						\{0\}	& \text{otherwise},
					\end{cases}
			\end{align*}
			Note that \(\lambda_k'\) is a constant in the definition of \(\psi_k\), and plays no r\^ole in the determination of the set \(\Psi(t, k)\). It does, however, influence the value of \(\ol\psi_k(t)\) at each \(t\in S\), and therefore, also the set-valued map \(t\mapsto \LQmuxset(t)\) defined in \eqref{e:LQmuxset def}.

			The numbers \(\bigl(\ol\psi_k(t)\bigr)_{k=1}^N\) defined in \eqref{e:olpsik def} admit the following concrete description:
			\begin{align*}
				\ol\psi_k(t)	& = \begin{cases}
					\lambda_k' + \frac{1}{2}\norm{\subadj{k}(t)}_{B_k R_k\inverse B_k\trnsp}^2	& \text{if }\norm{\subadj{k}(t)}_{B_k R_k\inverse B_k\trnsp}^2 \ge  2\lambda_k,\\
					\wt\lambda	& \text{otherwise},
				\end{cases}
				\\
					& \text{for all }t\in S.
			\end{align*}
			Finally, for each \(t\in S\),
			\begin{align*}
				& \LQmux(t)\in \\
				& \argmax_{k\in\{1, \ldots, N\}} \begin{cases}
						\lambda_k' + \frac{1}{2}\norm{\subadj{k}(t)}_{B_k R_k\inverse B_k\trnsp}^2	& \text{if }\norm{\subadj{k}(t)}_{B_k R_k\inverse B_k\trnsp}^2 \ge  2\lambda_k,\\
						\wt\lambda	& \text{otherwise}.
					\end{cases}
			\end{align*}
			Measurability of \(\LQmux\) follows from the assumption that \(\con\opt\) is measurable. The steps above immediately lead to the assertion.\hfill{}\(\square\)

%%%%%%%%%%%%%%%%%%%%%%%%%%%%%%%%%%%%%%%%%%%%%%%%%%%%%%%%%%%%%%%%%%%%%%%%%%%%%%%%
\section{Multiplexed sparse Mayer problem: proof}
\label{s:proof:Mayer}
%%%%%%%%%%%%%%%%%%%%%%%%%%%%%%%%%%%%%%%%%%%%%%%%%%%%%%%%%%%%%%%%%%%%%%%%%%%%%%%%
Recall that in \secref{s:Mayer} we defined \(\lambda_k' = \sum_{\ell\in\{1, \ldots, N\}\setmin\{k\}} \lambda_\ell\) for each \(k = 1, \ldots, N\), and \(\wt\lambda = \sum_{k=1}^N \lambda_k\).
\textbf{Proof of Theorem \ref{t:Mayer}}
	Theorem \ref{t:Clarke extended PMP} applies to \eqref{e:Mayer} because
	\begin{itemize}[label=\(\circ\), leftmargin=*]
		\item \(\admactMayer\) is a finite union of compact sets, and is therefore Borel measurable;
		\item the dynamics \(f\) is given by the linear control system \eqref{e:joint system} and is therefore smooth;
		\item the instantaneous cost function \(\Lambda(\xi, \mu) = - \sum_{k=1}^N \lambda_k \indic{\{0\}}(\mu^k)\) is independent of the space variable \(\xi\) and is lower semicontinuous in \(\mu\);
		\item the terminal cost \(\ell(\bar\xi, \hat\xi) = \frac{1}{2}\inprod{\hat\xi - \stfin}{\hat Q\bigl(\hat\xi - \stfin\bigr)}\) is smooth; and
		\item the boundary constraint set \(E = \{\stinit\}\times\R^d\) is closed.
	\end{itemize}

	Let the state-action trajectory \([\tinit, \tfin]\ni t\mapsto \bigl( \st\opt(t), \con\opt(t) \bigr)\) be a local minimizer of \eqref{e:Mayer}. For \(\eta\in\R\) we define the Hamiltonian
	\begin{equation}
		\label{e:Mayer:Hamiltonian}
		\begin{aligned}
				& (\R^{d_1}\times\cdots\times\R^{d_N})\times(\R^{d_1}\times\cdots\times\R^{d_N})\times\admactMayer\\
				& \quad \ni \left(\pmat{\subadj{1}\\ \smash[t]{\vdots} \\ \subadj{N}}, \pmat{\xi^1\\\smash[t]{\vdots}\\\xi^N}, \pmat{\mu^1\\\smash[t]{\vdots}\\\mu^N}\right) \teL (\adj, \xi, \mu) \mapsto\\
				& \quad \Ham^\eta(\adj, \xi, \mu) \Let \inprod{\adj}{A\xi + B\mu} + \eta\sum_{k=1}^N \lambda_k \indic{\{0\}}(\mu^k)\\
				& \quad = \sum_{k=1}^N \biggl( \inprod{\subadj{k}}{A_k \xi^k + B_k \mu^k} + \eta \lambda_k \indic{\{0\}}(\mu^k) \biggr) \in\R,
		\end{aligned}
	\end{equation}
	By Theorem \ref{t:Clarke extended PMP} there exists an absolutely continuous map \([\tinit, \tfin]\ni t\mapsto \adj(t)\in\R^d\), the adjoint trajectory, that solves
	\[
		- \dot \adj(t) = \partial_\xi \Ham^\eta\bigl(\adj(t), \cdot, \con\opt(t)\bigr)\bigl(\st\opt(t)\bigr) = A\trnsp \adj(t)\quad \text{for a.e.\ }t\in[\tinit, \tfin],
	\]
	or equivalently,
	\begin{equation}
		\label{e:Mayer:adjoints}
		-\tdsubadj{k}(t) = A_k\trnsp \subadj{k}(t)\quad\text{for a.e.\ }t\in[\tinit, \tfin],\quad k = 1, \ldots, N.
	\end{equation}
	The boundary conditions for the adjoint \(\adj\) are obtained from the transversality conditions \eqref{e:Clarke:transversality}, and for our problem \eqref{e:Mayer:adjoints} they are given by
	\[
		\bigl(\adj(\tinit), -\adj(\tfin)\bigr)\in\R^d\times\bigl\{\eta\hat Q\bigl(\st\opt(\tfin) - \stfin\bigr)\bigr\}.
	\]
	In other words, the differential equations for \((\subadj{k})_{k=1}^N\) have specific terminal boundary constraints and free initial conditions; therefore, they have to be solved in reverse time.

	Theorem \ref{t:Clarke extended PMP} admits only two values of \(\eta\). We claim that the case of \(\eta = 0\) does not arise in \eqref{e:Mayer}. Indeed, if \(\eta = 0\), then the final boundary condition of each \(\subadj{k}\) is \(0\), leading to \(\adj(\tfin) = 0\); since the linear adjoint differential equations in \eqref{e:Mayer:adjoints} have no forcing terms, the entire trajectory of \(\adj\) must vanish. This contradicts the non-triviality condition \eqref{e:Clarke:nontriviality} of Theorem \ref{t:Clarke extended PMP}. Therefore, \(\eta = 1\), and we write \(\Ham\) instead of \(\Ham^1\) henceforth.

	The optimal control actions as functions of time must satisfy the Hamiltonian maximum condition \eqref{e:Clarke:Hamiltonian max}: for a.e.\ \(t\in[\tinit, \tfin]\),
	\begin{align*}
		\con\opt(t) & \in \argmax_{\mu\in\admactMayer} \Ham\bigl(\adj(t), \st\opt(t), \mu\bigr)\\
					& = \argmax_{\mu\in\admactMayer} \biggl\{ \sum_{k=1}^N \inprod{B_k\trnsp \subadj{k}(t)}{\mu^k} + \lambda_k \indic{\{0\}}(\mu^k) \biggr\}.
	\end{align*}
	Proceeding as in the proof of Theorem \ref{t:reach} we see that
	\begin{align*}
		& \subcon{k}\opt(t) \in 
		\begin{dcases}
			\argmax_{v\in\subadmact{k}} \inprod{B_k\trnsp\subadj{k}(t)}{v} & \text{if} \max_{v\in\subadmact{k}} \inprod{B_k\trnsp\subadj{k}(t)}{v} > \lambda_k,\\
			\{0\}\cup\argmax_{v\in\subadmact{k}} \inprod{B_k\trnsp\subadj{k}(t)}{v} & \text{if} \max_{v\in\subadmact{k}} \inprod{B_k\trnsp\subadj{k}(t)}{v} = \lambda_k,\\
			\{0\} & \text{otherwise}.
		\end{dcases}
		\\
		& \qquad\text{for }k = 1, \ldots, N,
	\end{align*}
	and
	\[
		\Mayermux(t)\in\argmax_{k\in\{1, \ldots, N\}}\max_{v\in\subadmact{k}} \Bigl\{\inprod{B_k\trnsp \subadj{k}(t)}{v} + \lambda_k \indic{\{0\}}(v) + \lambda_k'\Bigr\}.
	\]
	Since \(\con\opt\) is measurable, so is \(\Mayermux\). The constancy of the Hamiltonian \eqref{e:Clarke:Hamiltonian constancy} gives the final assertion of the theorem, completing the proof.\hfill{}\(\square\)

\end{multicols}

\end{document}

%% file: intro.tex
%%%%%%%%%%%%%%%%%%%%%%%%%%%%%%%%%%%%%%%%%%%%%%%%%%%%%%%%%%%%%%%%%%%%%%%%%%%%%%%%
	\section{Introduction}
	\label{s:intro}
%%%%%%%%%%%%%%%%%%%%%%%%%%%%%%%%%%%%%%%%%%%%%%%%%%%%%%%%%%%%%%%%%%%%%%%%%%%%%%%%
		Let \(N\) be a positive integer, and consider the finite ensemble of linear time-invariant control systems given by
		\[
		\left\{
		\begin{aligned}
			& \tdsubst{k}(t) = A_k \subst{k}(t) + B_k \subcon{k}(t),\\
			& k = 1, \ldots, N.
		\end{aligned}
		\right.
		\]
		Suppose that this ensemble of linear systems are to be controlled in a way that at any time \(t\) only one of the controllers can be active while the rest of them are set to zero. In other words, there is a multiplexer or polling scheme that selects one system from the ensemble at each instant of time; the controller of the selected system may be non-zero while the other controllers are set to zero, resulting in the corresponding systems evolving in open-loop. In this article we are concerned with the design of sparse and optimal multiplexers.

		Multiplexing or polling arises naturally in situations where a central server must cater to a range of different tasks; if the server is incapable of parallel processing and is assigned the task of controlling \(N\) different systems, it must process them serially, leading to multiplexing. Alternatively, if the \(N\) controllers share a single communication channel that must be shared between them, the same problem of multiplexing the controllers arises. In this article we study three different control problems on \emph{sparse} and \emph{optimal multiplexing}. The first problem concerns ballistic reachability: Given a pair of initial and final states for each of the \(N\) linear systems and a time interval \(T > 0\), we synthesize, if possible, a multiplexing strategy for the controllers of the \(N\) systems such that at time \(T\) \emph{all} the \(N\) systems reach their designated final states. Of course, at any time during the interval \([0, T]\) only one of the controllers is permitted to be non-zero. In addition this ballistic reachability objective, we also \emph{simultaneously} stipulate that the control trajectories are as sparse as possible, i.e., for most of the time \emph{all} the control actions are set to \(0\).

		The second problem studied in this article is that of sparse and optimally multiplexed linear quadratic control. Given the \(N\) linear systems above, we synthesize optimally multiplexed controllers for the minimization of the sum of \(N\) quadratic performance objectives, one for each system. Sparsity of the control trajectories is enforced by \(\Lp{0}\)-regularization of the performance objectives, and optimality of the multiplexing is ensured by definition of the optimal control problem. The third problem is that of sparse Mayer problem. Recall that a Mayer-type optimal control problem involves only a terminal cost. While the standard Lagrange and Bolza forms of optimal control problems can easily recast into corresponding Mayer problems, the ones that we treat here are fine-tuned to playing surrogate to reachability problems. To wit, optimal control problems with a large penalty on the deviation of the terminal states away from a prespecified and desired terminal state provide good surrogates for ballistic reachability problems that are typically computationally difficult; however, such terminal cost Mayer problems are only approximations of ballistic reachability problems.

		The literature on the topic of optimal multiplexing appears be sparse. In the context of resource sharing over networks, \cite{ref:FarPriDel-14}, \cite{ref:FarPirDel-13} discuss multi-sensor fusion and scheduling algorithms in shared sensing actuator networks (SSAN's) while \cite{ref:VanSteBovHee-17} design switched PID controllers for SCARA robots to ensure reference tracking with bounded errors. The control and network co-design problem over interconnected systems via communication networks has been studied in \cite{ref:BarPhiZha-02}, \cite{ref:ZhaHri-06}, \cite{ref:SahBarMaj-15}, where the aim is to design network policies as well as control laws to compensate for packet losses and other network associated errors while maintaining stability of all subsystems. However, actuation or sensing resources are typically not assumed to be shared among the subsystems, and, in particular, there appears to have been no prior studies on optimal scheduling.

		Sparse controls is an emerging area in control theory, with a recent vigorous spurt of activity \cite{ref:FarLinJov-14}, \cite{ref:NagQueNes-16}, \cite{ref:JovLin-13}, \cite{ref:PolKhlSch-14}, \cite{ref:ChaNagQueRao-16}, \cite{ref:SriCha-16}. Sparsity in control is employed nowadays to improve the efficiency of electric engines by turning them off during periods of activity under a certain load threshold such as locomotive engines during \textsl{coasting}, in networked control systems where the time of use of communication channels need to be reduced, etc. Ideally, this problem consists of solving an \(\Lp{0}\)-optimal control problem, and the discrete-time versions of such problems are known to be computationally hard. The technical approach for continuous time sparse optimal control started, presumably motivated by the techniques in sparse signal processing, by considering \(\Lp{1}\)-approximations of \(\Lp{0}\)-optimal control problems, before it was noticed in \cite{ref:ChaNagQueRao-16} that the exact \(\Lp{0}\)-optimal control problem admits a crisp set of necessary conditions in the form of a non-smooth Pontryagin maximum principle. In this article we shall employ such non-smooth Pontryagin maximum principles to ensure sparsity of our controls. In the context of our multiplexed reachability problem, achieving maximal sparsity of the controls is our objective. In the context of our sparse linear quadratic and Mayer problems, we employ \(\Lp {0}\)-regularizations in the corresponding objective functions to ensure sparsity of the resulting control trajectories.

		The rest of the article unfolds as follows: \secref{s:prelims} contains a precise description of the systems under consideration. The sparse and optimal multplexing reachability, linear quadratic control, and Mayer problems are treated in \secref{s:reach}, \secref{s:LQ}, and \secref{s:Mayer}, respectively, with the proofs of the results presented in the Appendices. Numerical experiments are provided in \secref{s:examples}, where we demonstrate the effect of sparse and optimal multiplexing on an ensemble of two linear systems consisting of a controlled harmonic oscillator and a linearized inverted pendulum on a cart --- a \(2\)-dimensional and a \(4\)-dimensional system, respectively.

		The notations employed here are standard: Transposes of vectors and matrices are denoted by ``\({}\trnsp\)''. We let \(\inprod{v}{w} \Let v\trnsp w\) be the standard inner product on Euclidean spaces, \(\norm{\cdot}\) is the norm derived from this inner product by setting \(\norm{v} \Let \sqrt{\inprod{v}{v}}\), and let \(\norm{v}_M \Let \sqrt{\inprod{v}{Mv}}\) denote the \(M\)-weighted norm of \(v\) for a symmetric and non-negative definite matrix \(M\). Given two matrices \(M_1\) and \(M_2\), we let \(\blkdiag(M_1, M_2)\) denote the block-diagonal matrix \(\pmat{M_1 & 0\\ 0 & M_2}\). The indicator function \(\indic{S}\) of a set \(S\) is defined by \(\indic{S}(z) = 1\) if \(z\in S\) and \(0\) otherwise, and the set-valued map \(\sgn\) by
		\begin{equation}
		\label{e:sgn def}
			\R\ni y\mapsto \sgn(y) = \begin{dcases}
				\{1\}	& \text{if }y > 0,\\
				\{-1\}	& \text{if }y < 0,\\
				[-1, 1]	& \text{otherwise}.
			\end{dcases}
		\end{equation}
		We write \(\Leb\) for the Lebesgue measure on \(\R\).